\documentclass[10pt]{amsart}
\usepackage{amsfonts,epsfig,fancyhdr}
\usepackage{newlfont,amsfonts,amssymb,amsmath,amsthm,amsgen,amscd,enumerate}
\usepackage{multicol,picinpar}

\usepackage{euscript,color}

\def\Gbdle{\mathcal{G}}

\def\isom{\cong}
\def\galg{\mathfrak{g}}

\def\palg{\mathfrak{p}}
\def\halg{\mathfrak{h}}
\def\malg{\mathfrak{m}}

\def\som{\mathfrak{so}(\mathfrak{m})}

\def\dsum{\oplus}

\newcommand{\p}{\mathfrak{p}}
\newcommand{\Q}{\ensuremath{\mathbb{H}}}

\newcommand{\C}{\ensuremath{\mathbb{C}}}
\newcommand{\K}{\ensuremath{\mathbb{K}}}
\renewcommand{\O}{\ensuremath{\mathrm{O}}}
\newcommand{\R}{\ensuremath{\mathbb{R}}}
\renewcommand{\SS}{\ensuremath{\mathbb{S}}}
\newcommand{\vect}[1]{\mathbf{#1}}
\newcommand{\PSU}{\mathrm{PSU}}
\newcommand{\PSp}{\mathrm{PSp}}

\newcommand{\bmat}{\begin{pmatrix}}
\newcommand{\emat}{\end{pmatrix}}
\newcommand{\SSO}{\mathrm{SO}}
\newcommand{\CO}{\mathrm{CO}}
\newcommand{\SSp}{\mathrm{Sp}}
\newcommand{\SSU}{\mathrm{SU}}
\newcommand{\U}{\mathrm{U}}
\newcommand{\SSL}{\mathrm{SL}}
\newcommand{\GL}{\mathrm{GL}}
\newcommand{\1}{\mathbf{1}}
\newcommand{\ipq}{\1_{p,q}}

\newcommand{\e}{\mathrm{e}}
\renewcommand{\d}{{\mathrm d}}
\newcommand{\bcase}{\begin{case}}
\newcommand{\ecase}{\end{case}}

\newcommand{\bclaim}{\begin{claim}}
\newcommand{\eclaim}{\end{claim}}

\newcommand{\bstep}{\begin{step}}
\newcommand{\estep}{\end{step}}

\newcommand{\bhlem}{\begin{hlem}}
\newcommand{\ehlem}{\end{hlem}}

\newcommand{\Mat}{\mathrm{Mat}}

\newcommand{\bleer}{\begin{leer}}
\newcommand{\eleer}{\end{leer}}
\newcommand{\bde}{\begin{de}}
\newcommand{\ede}{\end{de}}

\newcommand{\ol}{\overline}

\newcommand{\mf}{\mathfrak}
\newcommand{\bs}{\begin{satz}}
\newcommand{\es}{\end{satz}}
\newcommand{\btheo}{\begin{theo}}
\newcommand{\etheo}{\end{theo}}
\newcommand{\bfolg}{\begin{folg}}
\newcommand{\efolg}{\end{folg}}
\newcommand{\blem}{\begin{lem}}
\newcommand{\elem}{\end{lem}}
\newcommand{\bnote}{\begin{note}}
\newcommand{\enote}{\end{note}}
\newcommand{\bprf}{\begin{proof}}
\newcommand{\eprf}{\end{proof}}
\newcommand{\bd}{\begin{displaymath}}
\newcommand{\ed}{\end{displaymath}}
\newcommand{\be}{\begin{eqnarray*}}
\newcommand{\ee}{\end{eqnarray*}}
\newcommand{\eeqa}{\end{eqnarray}}
\newcommand{\beqa}{\begin{eqnarray}}
\newcommand{\bi}{\begin{itemize}}
\newcommand{\ei}{\end{itemize}}
\newcommand{\bnum}{\begin{enumerate}}
\newcommand{\enum}{\end{enumerate}}
\renewcommand{\la}{\langle}
\renewcommand{\ra}{\rangle}

\newcommand{\beq}{\begin{equation}}
\newcommand{\eeq}{\end{equation}}
\newcommand{\einhalb}{\frac{1}{2}}
\newcommand{\rr}{\mathbb{R}}

\newcommand{\ccc}{\mathbb{C}}

\newcommand{\vf}{\varphi}
\newcommand{\earr}{\end{array}\]}
\newcommand{\barr}{\[\begin{array}}
\newcommand{\bvec}{\left(\begin{array}{c}}
\newcommand{\evec}{\end{array}\right)}

\newcommand{\sumi}{\sum_{i=1}^n}

\newcommand{\g}{\mathfrak{g}}
\newcommand{\h}{\mathfrak{h}}

\newcommand{\m}{\mathfrak{m}}
\newcommand{\n}{\mathfrak{n}}

\renewcommand{\sl}{\mathfrak{sl}}

\newcommand{\Hol}{\mathrm{Hol}}
\newcommand{\+}{\oplus}

\newcommand{\rrn}{\mathbb{R}^n}
\newcommand{\so}{\mathfrak{so}}
\renewcommand{\b}{\mathfrak{b}}

\newcommand{\gl}{\mathfrak{gl}}

\renewcommand{\sp}{\mathfrak{sp}}
\newcommand{\su}{\mathfrak{su}}
\renewcommand{\u}{\mathfrak{u}}
\newcommand{\w}{\omega}

\newcommand{\bbem}{\begin{bem}}
\newcommand{\ebem}{\end{bem}}
\newcommand{\bbez}{\begin{bez}}
\newcommand{\ebez}{\end{bez}}
\newcommand{\bbsp}{\begin{bsp}}
\newcommand{\ebsp}{\end{bsp}}

\newcommand{\trace}{\mathrm{tr }}

\newcommand{\wt}{\widetilde}

\newcommand{\tem}{\widetilde{M}}
\newcommand{\tg}{\widetilde{g}}

\newcommand{\Ad}{\mathrm{Ad}}
\newcommand{\ad}{\mathrm{ad}}
\newcommand{\G}{\mathrm{G}}

\theoremstyle{definition}
\newtheorem{de}{Definition}
\newtheorem{bem}{Remark}
\newtheorem{bez}{Notation}
\newtheorem{bsp}{Example}
\theoremstyle{plain}
\newtheorem{lem}{Lemma}
\newtheorem{satz}{Proposition}
\newtheorem{folg}{Corollary}
\newtheorem{theo}{Theorem}

\begin{document}
\bibliographystyle{abbrv}
\title[Conformal holonomy, symmetric spaces, and skew symmetric torsion]
{Conformal holonomy, symmetric spaces, and  skew symmetric torsion}

\dedicatory{Dedicated to Michael G. Eastwood on the occasion of his 60th birthday.}
\author{Jesse Alt}\address[Alt]{School of Mathematics, University of the Witwatersrand, Private Bag 3, Wits 2050, Johannesburg, South Africa}\email{jesse.alt@wits.ac.za}
\author{Antonio J. Di Scala}\address[Di~Scala]{Dipartimento di Scienze Matematiche,``G.L. Lagrange",  Politecnico di Torino, Corso Duca degli Abruzzi 24, 10129, Torino, Italia.}\email{antonio.discala@polito.it}
\author{Thomas Leistner}\address[Leistner]{School of Mathematical Sciences, University of Adelaide, SA 5005, Australia}\email{thomas.leistner@adelaide.edu.au}
\thanks{The authors acknowledge the support and hospitality by the International Erwin Schr\"odinger Institute in \mbox{Vienna}: Here the project was commenced during the workshop on ``Cartan Connections, Geometry of Homogeneous Spaces and Dynamics", and the results where first presented  at the workshop on ``The Interaction of Geometry and Representation Theory".
The last author was supported by the Australian Research Council through the fellowship FT110100429 and the grant DP120104582.}
\thanks{{\em {\bf Corresponding author:}} T.~Leistner, School of Mathematical Sciences, University of Adelaide, SA 5005, Australia, {\tt thomas.leistner@adelaide.edu.au}}

\begin{abstract}
We consider the question: Can the isotropy representation of an irreducible  pseudo-Riemannian symmetric space be realized as a conformal holonomy group? Using recent results by \v{C}ap, Gover and Hammerl, we study the representations of  $\SSO (2,1)$, $\PSU(2,1)$ and $\PSp(2,1)$ as isotropy groups of irreducible symmetric spaces of signature $(3,2)$, $(4,4)$ and $(6,8)$, respectively, describing the geometry induced by a conformal holonomy reduction to the corresponding subgroups. In the case of $\SSO (2,1)$,
we show that conformal manifolds with such a holonomy reduction are always locally conformally flat and hence this group cannot be a conformal holonomy group. This result completes the classification of irreducible conformal holonomy groups in Lorentzian signature. In the case of $\PSU(2,1)$, we show that conformal manifolds of signature $(3,3)$ with this holonomy reduction carry, on an open dense subset, a canonical nearly para-K\"ahler metric with positive Einstein constant. For $\PSp(2,1)$ we also show that there is an open dense subset endowed with a canonical Einstein metric in the conformal class. As a result, after restricting to an open dense subset, the conformal holonomy must be a proper subgroup of $\PSU(2,1)$ or of $\PSp(2,1)$, respectively. These are special cases of an interesting relationship between a class of special conformal holonomy groups, and non-integrable geometries with skew symmetric, parallel torsion, which we also explore. Finally, using a recent result of Graham and Willse we prove the following general non-existence result: for a real-analytic, odd-dimensional conformal manifold, the conformal holonomy group can never be given by the isotropy representation of an irreducible pseudo-Riemannian symmetric space unless the isotropy group is $\SSO^0(p+1,q+1)$.
\end{abstract}
\maketitle
\section{Introduction and statement of results}

A basic problem in differential geometry is to understand the holonomy representations of connections associated to geometric structures, e.g. to classify the holonomy representations which can be geometrically realized via a canonical connection. Given a linear connection on a principle fibre bundle, or equivalently, a covariant derivative $\nabla$ on a vector bundle  $\cal T$ over a manifold $M$, the holonomy group $\Hol_p(\cal T, \nabla)$ at a point $p\in M$ is defined as the group of parallel transports, with respect to $\nabla$, along loops in $M$ starting and ending at $p$.  The holonomy group is a subgroup of $\mathrm{GL}(\cal T_p)$ that inherits a Lie group structure from its connected component, the {\em restricted holonomy group}. The restricted holonomy is obtained by restricting the definition to contractible loops, and hence, both groups are the same for simply connected manifolds.

The most well-known case is of course the classification of Riemannian holonomies, accomplished in the 20th century as a result of work by many mathematicians, including most notably Berger \cite{berger55},  Alekseevsky, Calabi, Yau, Bryant and Joyce (see \cite{besse87}, \cite{bryant87}, \cite{joyce00} and references therein). For the classification of Lorentzian holonomy groups, see the survey \cite{galaev-leistner08} and references therein. Generalizations in various directions have been studied. For example, the classification of irreducible holonomy representations of torsion-free affine connections was completed by Merkulov and Schwachh\"ofer in \cite{merkulov-schwachhoefer99}\footnote{As the referee  pointed out to us, some representations missing from this classification were later noticed in \cite{bryant02}.}.

One possible generalization is to consider not a fixed pseudo-Riemannian metric, but only its conformal class, and study the holonomy associated to this geometric structure. When only a conformal equivalence class of pseudo-Riemannian metrics is fixed, the most natural connection to consider is not a principal connection but a {\em Cartan connection}. Given a Lie group $G$ with Lie algebra $\g$, a closed subgroup $P$ and a principle $P$-bundle $\cal G$, a Cartan connection $\w$ is a $P$-equivariant one-form on $\cal G$, that recovers fundamental vector fields and, in contrast to a principle fibre bundle connection, provides a global parallelism between $\cal G$ and $\g$ (see the definition in Section \ref{sec-background}). Hence, it does {\em not} define a horizontal distribution on $\cal G$, a fact which makes the notion of holonomy more involved. However, extending the bundle $\cal G$ to the $G$-bundle $\widehat{\Gbdle} := \Gbdle \times_P G$, the Cartan connection $\w$ induces a principle fiber bundle connection $\widehat{\w}$ on $\widehat{\cal G}$, and one can define the holonomy group of the Cartan connection $\w$ as the usual holonomy group of $\widehat{\w}$. There is a notion of holonomy for the Cartan connection that does not make use of this extension \cite[Section 5.4]{sharpe97}, but one can prove \cite[Proposition 1]{baum07} that both have the same connected component. This is one of the reasons why we restrict ourselves to the study of the restricted holonomy group of Cartan connections --- in the above sense as holonomy of $\widehat{\w}$. In the following, when we use the word holonomy we will always refer to the restricted holonomy group. The other reason is that our approach is based on the Lie algebra of the holonomy group, the {\em holonomy algebra}, which can only describe the connected component of the full holonomy. We also use the notion of {\em holonomy representation}, which refers to the restricted holonomy group, or its Lie algebra, and its representation on the fiber $\cal T_p$ of the vector bundle on which the covariant derivative $\nabla$ is defined.

The present work touches on the problem of classifying the representations which are realizable as the holonomy group of the canonical Cartan connection in conformal geometry, i.e., as ``conformal holonomy groups.'' This means we study the holonomy representations of the canonical (normal) Cartan connection induced by a conformal manifold of dimension at least $3$. (Definitions and some other relevant background material are reviewed in Sections 2.1--2.3.) In particular, note that this means that a conformal manifold $(M,[g])$ of signature $(p,q)$ has conformal holonomy $\mathrm{Hol}(M,[g])$ given as a subgroup of $\mathrm{O}(p+1,q+1)$, and thus the basic holonomy representation is on the space $\R^{p+1,q+1} \simeq \cal T_p$.

The study of conformal holonomy has attracted considerable interest in recent years. The first fundamental observation that was made is  that the conformal holonomy is contained in the stabilizer in $\O(p+1,q+1)$ of a line if and only if  there exists an  Einstein metric in the conformal class, where the Einstein metric might only be defined on an open dense subset \cite{leitner05,gover/nurowski04, leistner05a}. As a result, we know that if the conformal class $[g]$ contains an Einstein metric, then the conformal holonomy representation preserves a line in $\R^{p+1,q+1}$, so it does not act irreducibly. Another fundamental result discovered about conformal holonomy was an analog of the de Rham--Wu decomposition theorem, relating the decomposition of $\R^{p+1,q+1}$ into a direct sum of $\mathrm{Hol}(M,[g])$-invariant, non-degenerate subspaces to the existence of a metric in the conformal class $[g]$, again defined on a dense subset of $M$, which is locally a product of Einstein metrics (see \cite{armstrong07conf} and \cite{armstrong-leitner11} for Riemannian conformal structures and  \cite{leitner04} for arbitrary signature).  As a result, the case of conformal holonomies which act decomposably on $\R^{p+1,q+1}$ is fairly well understood. Our main focus in this work will be on the possible \emph{irreducible} representations on $\R^{p+1,q+1}$ which can be realized as conformal holonomy representations, which is the class of holonomies we can reasonably hope to make some progress toward classifying in general. A first step in this attempt was made in \cite{alt11}, where a classification is given under the additional assumption that the conformal holonomy acts transitively on the M\"obius sphere. The case of a degenerate subspace of dimension 2 was studied in \cite{leistner-nurowski12}, where it was shown that this corresponds to a pure radiation metric in the conformal class.

Note that in Riemannian signature irreducible representations play no significant role in any such classification, since the only connected subgroup of $\mathrm{O}(n+1,1)$ which acts irreducibly on $\R^{n+1,1}$ is the connected component $\SSO^0(n+1,1)$ (cf. \cite{discala-olmos01}, \cite{att05}, \cite{boubel-zeghib04}). For the classification of (non-irreducible) conformal holonomy representations in Riemannian signature, see  \cite{armstrong07conf,armstrong-leitner11}.

Turning next to Lorentzian signature, the following classification result was obtained by the second and third authors:

\btheo[{\cite[Corollary~1]{discala-leistner08}}]
\label{Lorentz signature class}
Let $H \subset \mathrm{O}(n,2)$ be a connected conformal holonomy group of an $n$-dimensional Lorentzian conformal manifold. If $H$ acts irreducibly on $\R^{n,2}$, then it must be one of the following:
$H = \SSO^0(n,2)$; $H=\SSU(m,1)$ for $n=2m$; or $H=\SSO^0(2,1)$ for $n=3$.
\etheo

Note that the last of these groups comes from the isotropy representation $$\Ad_{\SSL_3\R}: \SSO(2,1) \rightarrow \SSO(3,2)$$ of the (irreducible) pseudo-Riemannian symmetric space $\SSL_3\R/\SSO(2,1)$. The present work began with the observation that this representation can, in fact, be eliminated as a conformal holonomy group in Lorentzian signature. Indeed, this is our first main result:

\btheo \label{gensotheo}
For the irreducible pseudo-Riemannian symmetric space $\SSL_3\R/\SSO(2,1)$ of signature $(3,2)$, let $H \subset \SSO(3,2)$ be the image of $\SSO(2,1)$ under the isotropy representation. If a conformal manifold $(M,[g])$ has a conformal holonomy reduction to $H \subset \SSO(3,2)$, then $(M,[g])$ is locally conformally flat. In particular, its conformal holonomy group must be discrete, and thus the isotropy representation of $\SSL_3\R/\SSO(2,1)$ cannot be realized as a conformal holonomy group.
\etheo

As a corollary of Theorems \ref{Lorentz signature class} and \ref{gensotheo}, there are only two possible irreducible conformal holonomy representations in Lorentzian signature. The geometry in these two cases is well understood: The case $\SSO^0(n,2)$ corresponds to generic Lorentzian conformal manifolds; while $\SSU(m,1)$ corresponds, locally, to Lorentzian conformal manifolds which are the Fefferman space of some strongly pseudo-convex Cauchy-Riemann (CR) manifold of real dimension $(2m-1)$ (cf. \cite{fefferman76}, \cite{CapGover1}, \cite{CapGover2}).

The proof of Theorem \ref{gensotheo} can be seen as a basic application of the recent work of \v{C}ap, Gover and Hammerl in \cite{cgh11} which greatly clarifies the meaning of holonomy reduction for Cartan connections via the notion of ``curved orbit decompositions''. We will review this material in Section~\ref{sec-background}, along with other relevant background on Cartan geometries. Using the notion of ``curved orbit decomposition'' from \cite{cgh11}, a basic ingredient in the study of the geometry induced by the holonomy reduction of a Cartan geometry of type $(G,P)$ to some subgroup $H \subset G$ is an analysis of the $H$-orbits in the homogeneous model $G/P$. In particular, for conformal Cartan geometry the homogeneous model $G/P$ is a double covering of the M\"obius sphere $\SS^{p,q}$, where the latter is viewed as the projectivization of the null cone in $\R^{p+1,q+1}$. Thus we are led to study the action of the isotropy subgroup of the symmetric spaces $\SSL_n\R/\SSO(p,q)$ on the space of null lines in $$T_o(\SSL_n\R/\SSO(p,q)) \cong \sl_n\R/\so(p,q) \cong \m = \{ X \in \sl_n\R : \la Xu,v\ra=\la u,Xv\ra \}.$$ This is the subject of Section \ref{sec-orbit}, where we also look at the closely related symmetric spaces $\SSL_n \C / \SSU(p,q)$ and $\SSL_n \mathbb{H} / \SSp(p,q)$. In these cases, we show that open orbits in the relevant M\"obius sphere can occur only for $n=3$ (for example,  for $\SSL_n\R/\SSO(p,q)$ this is so because $\frac{n(n-1)}{2}=\mathrm{dim}(\SSO(p,q))$ is smaller than $\frac{(n-1)(n+2)}{2} - 2$ if $n>3$, which is the dimension of the relevant M\"obius sphere, cf. Figure \ref{dimension table}), and for $n=3$ the union of open orbits is dense. In terms of the geometry induced by a holonomy reduction, it is thus of primary interest to look at the isotropy representations for $n=3$.

The isotropy representations of $\mathrm{SL}_3 \C /\SSU(2,1)$ and $\mathrm{SL}_3 \mathbb{H} /\SSp(2,1)$ give irreducible representations $\PSU(2,1) \subset \SSO(4,4)$ and $\PSp(2,1) \subset \SSO(6,8)$, respectively. In Section 4, we apply the results of \cite{cgh11} and the analysis of Section \ref{sec-orbit} to the question of the realizability of these conformal holonomy representations. After proving Theorem \ref{gensotheo}, we obtain the following result for the case of $\PSU(2,1) \subset \SSO(4,4)$:

\btheo\label{sutheo}
If $(M,[g])$ is a conformal manifold of signature $(3,3)$ with conformal holonomy $\mathrm{Hol}(M,[g]) \subseteq \PSU(2,1) \subset \SSO(4,4)$, then, on an open dense subset $M_0 \subset M$, there exists a canonical metric with nearly para-K\"ahler structure. In particular,
\bnum
\item on $M_0$, there is an Einstein metric $g_0 \in [g]$ with positive Einstein constant, and
\item  the conformal holonomy $\mathrm{Hol}(M_0, [g])$ preserves a time-like vector in $\rr^{4,4}$, and hence is {\em properly} contained in $\PSU(2,1)$.
\enum
\etheo

The existence of the nearly para-K\"ahler structure is perhaps just as interesting as the fact that the conformal holonomy group $\PSU(2,1) \subset \SSO(4,4)$ can be excluded, at least if one restricts to an open dense submanifold. It is induced by a certain metric affine connection with skew-symmetric and parallel torsion, which comes quite naturally from analyzing the homogeneous orbits and the resulting holonomy reduction of the normal conformal Cartan connection. In fact, the natural properties of this affine connection imply that its geometry is in a certain sense quite ``close'' to that of the naturally reductive homogeneous geometry which the Cartan geometry inducing it is modeled on, specifically that their Ricci tensors are equal. This idea also works in the case of $\PSp(2,1) \subset \SSO(6,8)$, leading to the:

\btheo\label{sptheo}
If $(M,[g])$ is a conformal manifold of signature $(5,7)$ with conformal holonomy $\mathrm{Hol}(M,[g]) \subseteq \PSp(2,1) \subset \SSO(6,8)$, then there is an open dense subset $M_0 \subset M$ and a canonical Einstein metric $g_0 \in [g |_{M_0}]$. In particular, the conformal holonomy $\mathrm{Hol}(M_0,[g])$ preserves a line in $\R^{6,8}$ and hence is a \emph{proper} subgroup of $\PSp(2,1)$.
\etheo

These two examples are an instance of the more general principle established in Theorem \ref {Cartan invariant conformal holonomy} of Section \ref{sec-background}, and indicate an interesting connection between special conformal holonomy and metric connections with torsion (cf. \cite{agricola-srni} for a survey of the latter). The induced metric connections with skew-symmetric, parallel torsion given in both cases, are a rather intriguing twist on the idea expressed by \'E. Cartan: ``Given a manifold embedded in affine (or projective or conformal etc.) space, attribute to this manifold the affine (or projective or conformal etc.) connection that reflects in the simplest possible way the relations of this manifold with the ambient space'' (\cite{cartan24}, quoted from \cite{agricola-srni}). The twist being that in the cases we consider, the conformal Cartan geometry with special holonomy plays the role of ``ambient space,'' while the torsion of the distinguished metric affine connection enters naturally as a reflection of how the nearly para-K\"ahler metric, etc., lies in the conformal class having special holonomy.

Note that, under the assumption that the manifold and the conformal structure are real analytic, i.e., that there is an analytic metric in the conformal class, we obtain the result of Theorems~\ref{sutheo}~and~\ref{sptheo} globally, that is, without restricting to an open dense subset\footnote{We should mention explicitly that we have not attempted to address the much more difficult question of whether non-analytic conformal manifolds can be found with the full holonomies $\PSU(2,1)$ or $\PSp(2,1)$.}. This is based on the well-known fact \cite[Section II.10]{ko-no63} that the holonomy algebra at $p\in M$ of a linear  connection on a principle fiber bundle over $M$ is equal to the {\em infinitesimal holonomy algebra}, which is defined as  the span of all derivatives of the curvature {\em at the point $p$}, provided that the bundle and the connection are real analytic. Since we have defined the conformal holonomy as the holonomy of the associated principle fiber bundle  connection, this yields
\bfolg
Let $(M,[g])$ be a real analytic conformal manifold of signature $(3,3)$, respectively $(5,7)$. Then its conformal holonomy representation cannot be the isotropy representation of the irreducible pseudo-Riemannian symmetric symmetric space $\mathrm{SL}_3 \C /\SSU(2,1)$, respectively $\mathrm{SL}_3 \mathbb{H} /\SSp(2,1)$.
\efolg

From the limited evidence available, a natural question to ask is whether it is ever possible to realize the full isotropy subgroup of an irreducible pseudo-Riemannian symmetric space as a conformal holonomy representation. A first partial answer can be obtained as a corollary to the following statement proved in  \cite{CapGover2} and \cite{leitner07}: If the conformal holonomy of a conformal manifold of signature $(2p+1,2q+1)$ is contained in $\U(p+1,q+1)$, then it is already contained in  $\SSU(p+1,q+1)$. On the other hand, the fact that irreducible symmetric spaces cannot be Ricci-flat without being flat (cf. \cite{alekseevsky10} or Proposition \ref{Ricci nonflat} below) prevents them from having holonomy contained in the special unitary group. Hence we have
\bfolg
Let $H\subset \U(p+1,q+1)$ act irreducibly and be given as the isotropy representation of a non-flat pseudo-K\"ahler symmetric space. Then $H$ cannot be a conformal holonomy group.
\efolg
As a further partial answer to the above question, in Section \ref{fg} we prove:
\btheo \label{ambient holonomy extension}
Let $(M,[g])$ be a real-analytic conformal manifold of signature $(p,q)$ and odd dimension $n=p+q \geq 3$. If $H \subset \O(p+1,q+1)$ acts irreducibly on $\rr^{p+1,q+1} $ and is defined as the identity component of the stabilizer in $\O(p+1,q+1)$ of some tensor, then $H$ cannot be equal to the conformal holonomy group of $(M,[g])$ unless $H=\SSO^0(p+1,q+1)$ or $H=\mathrm{G}_{2(2)}\subset \SSO(4,3)$, where $\mathrm{G}_{2(2)}$ is the split real form of the simple Lie group $\mathrm{G}_2$.
\etheo
The proof relies on the Fefferman-Graham ambient metric construction of conformal geometry \cite{fefferman/graham85,fefferman-graham07} and a related result by Graham and Willse \cite{graham-willse11} which, in odd dimensions and with the assumption of real analyticity,  guarantee that  a unique  Ricci flat ambient metric extists and that parallel tractors extend to parallel ambient tensors. Then the theorem essentially follows from Berger's list of non-symmetric, irreducible, pseudo-Riemannian holonomy groups, and the observation that Ricci-flat manifolds cannot have the holonomy of a pseudo-Riemannian irreducible symmetric space unless $H=\SSO^0(p+1,q+1)$. Note that candidates for conformal structures with conformal holonomy equal to $\mathrm{G}_{2(2)}$ were discovered by Nurowski \cite{nurowski04}. Their ambient metric was first studied in \cite{nurowski07}, and in \cite{leistner-nurowski11} and \cite{graham-willse11} conditions on the conformal structures were given for which the ambient metric has holonomy equal to $\G_{2(2)}$. Note, however, that it has not yet been verified for any of those examples if the \emph{conformal holonomy} equals $\G_{2(2)}$ (in general the conformal holonomy is contained in the holonomy of the ambient metric, but {\em a priori} this containment could be proper), but we believe that this is only a matter of computing sufficiently many derivatives of the tractor curvature.\footnote{This equality was established recently by \v{C}ap, Gover, Graham, and Hammerl in their work in progress \cite{cggh12}. They prove a general result relating the ambient and conformal holonomy which also implies our  Theorem 5.}

Returning to the question of realizability of isotropy representations of symmetric spaces as conformal holonomy groups, we recall that isotropy groups of irreducible pseudo-Riemannian symmetric spaces are given as stabilizer of their curvature tensor at a point (cf. Proposition \ref{stabilizer} below) and obtain:

\bfolg
\label{Thomas conjecture real analytic odd}
Let $(M,[g])$ be a real-analytic conformal manifold of signature $(p,q)$ and odd dimension $n=p+q \geq 3$. If $H \subset \SSO^0(p+1,q+1)$ is a connected, irreducibly acting proper subgroup given by the isotropy representation of an irreducible pseudo-Riemannian symmetric space and $\mathrm{Hol}(M,[g]) \subseteq H$, then $\mathrm{Hol}(M,[g])$ is a \emph{proper} subgroup of $H$.
\efolg

Theorems \ref{gensotheo}, \ref{sutheo} and \ref{sptheo}, and Corollary \ref{Thomas conjecture real analytic odd}, lead us to surmise that a non-existence result could likely be true in general, but proving or disproving such a result is beyond the scope of the present work. It is notable that this would be a remarkable contrast to the situation in Riemannian holonomy, where the full isotropy subgroups are realizable as holonomy representations of locally symmetric spaces, but are not strictly excluded. Indeed, one has the impression that irreducible conformal holonomy representations are much more ``scarce'' than in the Riemannian case, an impression supported by other examples of representations which are not realizable as conformal holonomy groups. For example, by the result in \cite{leitner07} and \cite{CapGover2} mentioned above, the standard representation of the indefinite unitary group, $\mathrm{U}(p+1,q+1) \subset \SSO(2p+2,2q+2)$, cannot be realized as a conformal holonomy group: there is no ``conformal analog'' of K\"ahler manifolds as distinct from Calabi-Yau manifolds. Moreover, to our knowledge, no example of a conformal holonomy group that is {\em not} a pseudo-Riemannian holonomy group has been found.  However, apart from the result in \cite[Theorem 3.2]{leistner05a} (see a related result in \cite{armstrong-leistner07}) that the conformal holonomy of a conformal $C$-space is a Berger algebra, there is as yet no known analog of the Berger criteria for Riemannian holonomy which could give effective (algebraic) restrictions on the possible irreducible conformal holonomy representations. One might hope that the methods needed to prove a conjecture excluding isotropy representations, would also yield some insight into the appropriate conformal Berger criteria.
\subsection*{Acknowledgements} The authors would like to thank the referees for their valuable observations and comments.

\section{Cartan geometry and pseudo-Riemannian symmetric spaces}
\label{sec-background}
\subsection{Cartan connections} In this section, we review some relevant facts about holonomy reductions for general Cartan geometries (Section \ref{sec-background2}), the Cartan geometry corresponding to conformal structures (Section \ref{sec-background3}), and Cartan geometries of reductive type (Section \ref{sec-background4}) and apply this to the holonomy reductions to isotropy groups of irreducible symmetric spaces (Section \ref{sec-background5}). Recall that a Cartan geometry $(\pi: \Gbdle \rightarrow M,\omega)$ of some type $(G,P)$ is given, for $G$ a Lie group and $P$ a closed subgroup, by a $P$-principal bundle $\pi: \Gbdle \rightarrow M$ and a one-form with values in the Lie algebra of $G$, $\omega \in \Omega^1(\Gbdle;\galg)$, satisfying the axioms of a \emph{Cartan connection}, i.e. $\omega$:
\begin{itemize}
  \item trivializes the tangent bundle of the total space: $\omega_u: T_u\Gbdle \stackrel{\simeq}{\rightarrow} \mathfrak{g}$ for all $u \in \Gbdle$;
  \item is $P$-equivariant: $R_p^*\omega = \mathrm{Ad}(p^{-1}) \circ \omega$;
  \item recovers fundamental vector fields: $\omega(\widetilde{X}) = X$ for all $X \in \palg$, where \[\widetilde{X}(u) := \frac{d}{dt}\vert_{t=0}(u.\mathrm{exp}(tX))\] denotes the fundamental vector field.
\end{itemize}
For background on the theory of Cartan geometries, the reader is referred to the books \cite{sharpe97} and \cite{cap-slovak-book09}. A viewpoint which is often useful to take is of the Cartan geometry $(\Gbdle,\omega)$ being a ``curved version'' of the homogeneous model geometry $G \rightarrow G/P$ (in which $G$ is the automorphism group determining a ``geometry'' on $G/P$), with the Cartan connection $\omega$ playing the role of the Maurer-Cartan form $\omega_G: T_gG \ni X \mapsto (L_{g^{-1}})_*(X) \in \galg$.

The \emph{curvature $2$-form} of $(\Gbdle,\omega)$ is defined to be $\Omega^{\omega} := d\omega + \einhalb[\omega,\omega] \in \Omega^2(\Gbdle;\galg)$, and it is a well-known fact that $\Omega^{\omega}$ vanishes identically if and only if $(\Gbdle,\omega)$ is locally isomorphic to the homogeneous model geometry $(G \rightarrow G/P,\omega_G)$ (for a proof, see Chapter 5 of \cite{sharpe97}). Such a Cartan geometry (with vanishing curvature form) is called \emph{flat}. It is often useful to translate the curvature form into a curvature function $\kappa^{\omega} \in C^{\infty}(\Gbdle;\Lambda^2(\g/\p)^* \otimes \g)$, by defining, for $u \in \Gbdle$ and $X,Y \in \g$, $\kappa^{\omega}(u)(X,Y) := \Omega^{\omega}(\omega_u^{-1}(X),\omega_u^{-1}(Y))$, and noting that $\kappa^{\omega}(u)(X,Y)=0$ whenever $X \in \p$, since $\Omega^{\omega}$ is horizontal.

We also mention here that composing the curvature tensor $\Omega^{\omega}$ with the natural projection $\galg \rightarrow \galg/\palg$ defines the \emph{torsion} $\Theta^{\omega} \in \Omega^2(\Gbdle;\galg/\palg)$ of the Cartan connection, and a Cartan geometry is called \emph{torsion-free} if its torsion vanishes identically, i.e., if and only if $\Omega^{\omega} \in \Omega^2(\Gbdle;\palg)$, or equivalently, if  $\kappa^{\omega} \in C^{\infty}(\Gbdle;\Lambda^2(\g/\p)^* \otimes \p)$. Obviously, any flat Cartan geometry is torsion-free, but there are also many important examples of torsion-free Cartan geometries which need not be flat, for instance the canonical (normal) Cartan geometry of a conformal manifold (cf. Section \ref{sec-background3}) is always torsion-free.

\subsection{Holonomy reduction and curved orbit decompositions for Cartan geometries}
\label{sec-background2}
There is a well-defined notion of \emph{holonomy} of a Cartan geometry, $\mathrm{Hol}(\Gbdle,\omega)$, determined up to conjugation as a subgroup of $G$. This is obtained by taking the extension $\widehat{\Gbdle} := \Gbdle \times_P G$ to a $G$-principal bundle, noticing that there is a unique principal connection $\widehat{\omega} \in \Omega(\widehat{\Gbdle};\galg)$ which pulls back to $\omega$ under the natural inclusion $\iota: \Gbdle \hookrightarrow \widehat{\Gbdle}$, and defining 
$$\mathrm{Hol}(\Gbdle,\omega) := \mathrm{Hol}(\widehat{\Gbdle},\widehat{\omega}),$$ 
where $\mathrm{Hol}(\widehat{\Gbdle},\widehat{\omega})$ is defined in the usual way for principal bundle connections, via horizontal path-lifting.

As usual with principal bundle connections, useful tools for studying the holonomy of a Cartan connection are obtained from associated bundles. In particular, a linear representation $\rho: G \rightarrow \GL(\mathbb{W})$ determines the associated vector bundle $W(M) := \widehat{\Gbdle} \times_{\rho} \mathbb{W}$ which inherits in the usual way a covariant connection $\nabla^W$ canonically induced by $\widehat{\omega}$ (and hence by $\omega$). $W(M)$ and $\nabla^W$ are called \emph{tractor bundle} and \emph{tractor connection} of the Cartan geometry. Note that $W(M)$ can also be considered as an associated vector bundle to the $P$-bundle $\Gbdle \rightarrow M$, by restricting the representation $\rho$ to $P$, $W(M) = \Gbdle \times_{\rho(P)} \mathbb{W}$.

Suppose a subgroup $H  \subset G$ is given as the automorphism subgroup of some element $\alpha \in \mathbb{W}$, $H = \mathrm{Aut}(\alpha) := \{ g \in G : \rho(g).\alpha = \alpha \}$. Then from standard facts about principal bundle connections, we know that $\mathrm{Hol}(\Gbdle,\omega) \subseteq H$ if and only if there exists some $\nabla^W$-parallel section (tractor) $s \in \Gamma(W(M))$ \emph{of type $\alpha$} (that is, such that $s(x) = [(u,\alpha)]$ for any $x \in M$ and some extended frame $u \in \widehat{\Gbdle}_x$; identifying $s \in \Gamma(W(M))$ with a $G$-equivariant function $s: \widehat{\Gbdle} \rightarrow \mathbb{W}$, this can be written $\alpha \in s(\widehat{\Gbdle}_x)$). Evidently, from equivariance of the constructions involved, we could just as well replace $\alpha$ by any other $\alpha' \in \mathcal{O} := \rho(G).\alpha \subset \mathbb{W}$. Then the condition that $s$ has type $\alpha$ can be replaced by the $G$-invariant condition $s(\widehat{\Gbdle}_x) = \mathcal{O}$, and we call the $G$-orbit $\mathcal{O}$ the \emph{$G$-type of $s$ at $x \in M$}. If $s$ is $\nabla^W$-parallel, then its $G$-type is evidently constant over $x \in M$ if $M$ is connected.

While these facts are more or less standard from the holonomy theory for principal connections, the problem is how to relate the holonomy reduction for $(\widehat{\Gbdle},\widehat{\omega})$ back to some construction for the Cartan geometry $(\Gbdle,\omega)$. A solution has recently been given by \cite{cgh11} via the concept of the \emph{curved orbit decomposition} determined by a parallel section $s \in \Gamma(W(M))$, cf. \cite[Section 2.4]{cgh11}. The key observation there is that we have a further (point-wise) invariant associated to $s$, which detects the fact that $(\widehat{\Gbdle},\widehat{\omega})$ was induced from $(\Gbdle,\omega)$: For $x \in M$, the \emph{$P$-type of $s$ at $x$} is defined to be the $P$-orbit $\overline{\alpha} := \rho(P).\alpha \subset \mathbb{W}$ such that $s(\Gbdle_x) = \overline{\alpha}$. By equivariance, this is well-defined for a fixed $x \in M$, but note that in contrast to the $G$-type, the $P$-type of a $\nabla^W$-parallel $s \in \Gamma(W(M))$ may change over $x \in M$. So one gets a decomposition of the base space $M$ (the \emph{curved orbit decomposition}) as:
\begin{align}
M = \bigsqcup_{\overline{\alpha} \in P \backslash \mathcal{O}} M_{\overline{\alpha}}, \label{curved orbit decomposition}
\end{align}
where $M_{\overline{\alpha}} := \{ x \in M : s(\Gbdle_x) = \overline{\alpha} \}$. The different possible $P$-types are indexed by the orbit space $P \backslash \mathcal{O} \cong P \backslash G/H$, but one has an isomorphism $P \backslash G/H \isom H \backslash G/P$ (cf. (3) of \cite{cgh11}), so each possible $P$-type corresponds precisely to an $H$-orbit in the homogeneous model $G/P$. In fact, the decomposition (\ref{curved orbit decomposition}) just gives, in the case of the homogeneous Cartan geometry $(G \rightarrow G/P,\omega_G)$ and its holonomy reduction to an automorphism subgroup $H = \mathrm{Aut}(\alpha) \subset G$, the decomposition into $H$-orbits, whence the name ``curved orbit decomposition.'' Now the main result of \cite{cgh11} which we will make use of can be stated as:

\btheo
[\v{C}ap, Gover \& Hammerl {\cite[Theorem 2.6]{cgh11}}]
\label{cgh theorem}
 Let $(\pi: \Gbdle \rightarrow M,\omega)$ be a Cartan geometry of type $(G,P)$, and suppose $\mathrm{Hol}(\Gbdle,\omega) \subseteq H \subset G$ for $H = \mathrm{Aut}(\alpha)$, some $G$-module $\mathbb{W}$ and $\alpha \in \mathbb{W}$. Then $M$ decomposes according to $P$-types into a disjoint union of initial submanifolds $M_{\overline{\alpha}}$. For any points $x \in M_{\overline{\alpha}}$ and $y \in \overline{\alpha} \subset G/P$, there exist neighborhoods $U_x \subset M$ of $x$ and $V_y \subset G/P$ of $y$ and a diffeomorphism $\varphi_{\overline{\alpha}}: V_y \rightarrow U_x$ such that $\varphi_{\overline{\alpha}}(V_y \cap \overline{\alpha}) = U_x \cap M_{\overline{\alpha}}$ and such that the diagram
\be
G/P\supset V_y&\stackrel{\varphi_{\overline{\alpha}}}{\to}& U_x\subset M
\\
\downarrow\ &&\ \downarrow
\\
H\backslash G/P&\simeq & P\backslash G/H,
\ee
in which the downward arrows assign to every point in the neighbourhoods its orbit type,
commutes.

Moreover, each curved orbit $M_{\overline{\alpha}}$ carries a naturally induced Cartan geometry $(\pi: \mathcal{H}_{\overline{\alpha}} \rightarrow M_{\overline{\alpha}}, \omega_{\overline{\alpha}})$ of type $(H,P_{\overline{\alpha}})$ for $H/P_{\overline{\alpha}} \isom \overline{\alpha} \subset G/P$. This geometry reduces $(\Gbdle,\omega)$, in particular its curvature $\Omega^{\omega_{\overline{\alpha}}}$ is given by restricting $\Omega^{\omega}$ to a sub-bundle $\mathcal{H}_{\overline{\alpha}}$ of $\Gbdle$ over $M_{\overline{\alpha}}$, and it is torsion-free whenever $(\Gbdle,\omega)$ is.
\etheo
\label{sec-background2}

\subsection{The normal Cartan geometry of a conformal manifold} 
\label{sec-background3} The conformal holonomy group which is the subject of this article, is defined, using the general approach reviewed in Section \ref{sec-background2}, from the canonical normal Cartan geometry associated to a conformal manifold. We now review briefly the basic facts about this Cartan geometry. For more background and proofs of many of the facts cited here, see for example Chapter 8 of \cite{sharpe97} or Chapter 1.6 of \cite{cap-slovak-book09}.

Let $M$ be a smooth manifold of dimension $n$ and $(p,q)$ some non-negative integers with $p+q=n$. Let $G := \O(p+1,q+1)$, $e_+ \in \R^{p+1,q+1}$ some non-zero null vector, and define $P \subset G$ to be the (closed, parabolic) subgroup which preserves the null ray $\R_+ e_+ \subset \R^{p+1,q+1}$ under the standard representation of $G$. Then Cartan geometries of type $(G,P)$ over $M$ correspond to conformal structures of signature $(p,q)$ on $M$.

One direction in this correspondence is not very difficult to establish, and we review it now briefly for future reference: A Cartan geometry $(\Gbdle \rightarrow M,\omega)$ of type $(G,P)$, as defined above, induces a conformal equivalence class $[g]$ of metrics of signature $(p,q)$ on $M$. To see this fact, consider the homomorphism $\overline{\mathrm{Ad}}: P \rightarrow \mathrm{GL}(\galg/\palg)$, given by the induced adjoint action on the quotient $\galg/\palg$, and the normal subgroup $P_+ := \mathrm{Ker}(\overline{\mathrm{Ad}}) \lhd P$. Then one can verify that $P/P_+ \cong \CO(p,q)$. Explicitly, if we fix another null vector $e_- \in \R^{p+1,q+1}$ which is dual to $e_+$, i.e. such that $\la e_+,e_-\ra = 1$, and let $P_0 \subset P$ be the subgroup which also preserves the null ray $\R_+ e_- \subset \R^{p+1,q+1}$, then calculating in a basis $\{e_+,e_1,\ldots,e_n,e_-\}$ with $\{e_1,\ldots,e_n\}$ an orthonormal sub-basis of $(\R e_+ \oplus \R e_-)^{\perp} \cong \R^{p,q}$ shows that $P_0 \cong \CO(p,q)$ and $P/P_+ \cong P_0$ (cf. e.g. Section 1.6.3 of \cite{cap-slovak-book09}).

Now the above facts can be used to define a conformal structure on $M$, using the following fact about the Cartan geometry $(\pi: \Gbdle \rightarrow M,\omega)$ (which is valid for arbitrary Cartan geometries, cf. 5.3 of \cite{sharpe97}): We have an isomorphism $TM \cong \Gbdle \times_{(\overline{\mathrm{Ad}},P)} \galg/\palg$; explicitly, for a point $x \in M$, any point $u \in \Gbdle_x$ determines an isomorphism
\begin{align}
\psi_u: T_xM \rightarrow \galg/\palg \label{Cartan frames}
\end{align}
by mapping $\psi_u: X \mapsto \omega(\widehat{X}) + \palg$,where $\widehat{X} \in T_u \Gbdle$ is any tangent vector which projects to $X$ via $\pi_*$. Moreover, the isomorphisms (\ref{Cartan frames}) satisfy the equivariance property, $\psi_{u.p} = \overline{\mathrm{Ad}}(p^{-1}) \circ \psi_u$, for any $p \in P$. In particular, since $\overline{\mathrm{Ad}}(P) \cong \CO(p,q) \subset \mathrm{GL}(\galg/\palg)$, the Cartan geometry $(\Gbdle \rightarrow M,\omega)$ determines a reduction of the structure group of $TM$ to the conformal group of signature $(p,q)$ and therefore a conformal class of metrics on $M$.

For applications, we often will also want to know how a choice of metric in the induced conformal class can be specified, and how to evaluate this metric on tangent vectors in terms of the Cartan connection. This is done via a choice of the dual null vector $e_- \in \R^{p+1,q+1}$ which was used to establish the isomorphism $P/P_+ \cong \CO(p,q)$ above. Namely, the calculation in the basis $\{ e_+,e_1,\ldots,e_n,e_- \}$ also shows that the subalgebra $\palg$ has an $\mathrm{Ad}(P_0)$-invariant complement $\p_- \subset \galg = \mathfrak{so}(p+1,q+1)$ and this identifies $\p_- \cong \g/\p \cong (\R e_+ \oplus \R e_-)^{\perp}$ as $P_0$-modules (the subalgebra $\p_- \subset \g$ can be defined by letting $\widehat{\p} := \mathfrak{stab}(\R e_-)$ and $\p_- := \mathrm{Ker}(\ol{\mathrm{ad}}: \widehat{\p} \rightarrow \mathfrak{gl}(\g/\widehat{\p}))$). In particular, the unique $P_0$-invariant conformal class of signature $(p,q)$ metrics on $\p_- \cong \galg/\palg$ is given by taking the conformal class of the isometric image of the natural metric on $(\R e_+ \oplus \R e_-)^{\perp} \cong \R^{p,q}$. Moreover, any choice of dual null vector $e_-$ determines a subgroup of $P$ which is isomorphic to $\O(p,q)$, namely the subgroup which fixes the vectors $e_+$ and $e_-$ (and not just the rays they determine). Now a metric in the conformal class defined by $(\Gbdle \rightarrow M,\omega)$ is easily seen to be determined by choosing a dual null vector $e_-$ and a reduction of $\Gbdle \rightarrow M$ to this copy of the pseudo-orthogonal group $\O(p,q)$ as a subgroup of $P$. To evaluate the resulting metric on tangent vectors $X, Y \in T_xM$, we use the isomorphisms (\ref{Cartan frames}), restricted to the reduced frames, and the corresponding $\O(p,q)$-invariant metric on $\galg/\palg$.

The surprising fact about conformal structures (which is also much more difficult to prove) is that the above construction can be ``reversed'', at least for $n \geq 3$: Given a conformal manifold $(M,[g])$ of signature $(p,q)$, with $p+q \geq 3$, there always exists a Cartan geometry $(\Gbdle \rightarrow M,\omega)$ of type $(G,P)$ which induces the initial conformal structure $[g]$. Moreover, this Cartan geometry can be uniquely determined, up to isomorphism, by a \emph{normalization} condition on the curvature of $\omega$: For $\g = \p_- \oplus \p_0 \oplus \p_+$ the decomposition as above, it is a general fact (which is not difficult to verify directly) that the subalgebras $\p_-$ and $\p_+$ admit bases $\{X_1,\ldots,X_n\}$ and $\{Z^1,\ldots,Z^n\}$, respectively, which are \emph{dual} with respect to the Killing form $K_{\g}$, i.e. they satisfy $K_{\g}(X_i,Z^j) = \delta_{ij}$; then normality of $\omega$ means that the composition $(\partial^* \circ \kappa^{\omega}) \in C^{\infty}(\mathcal{G};(\g/\p)^* \otimes \g)$ vanishes identically, i.e. that
\begin{align}
(\partial^* \circ \kappa^{\omega})(u)(X) &:= \sum_{i=1}^n [\kappa^{\omega}(u)(X,X_i),Z^i]_{\g} = 0 \label{normality}
\end{align}
for all $u \in \mathcal{G}$ and all $X \in \g$. The identity (\ref{normality}) is $P$-invariant in $u$ and independent of the choice of $K_{\g}$-dual bases of $\p_-$ and $\p_+$. We call the Cartan geometry of type $(G,P)$ which induces the conformal structure $(M,[g])$ and satisfies (\ref{normality}) the \emph{normal (or canonical) conformal Cartan geometry/connection} of $(M,[g])$. The theorem stating that the normal conformal Cartan connection exists and is unique up to isomorphism is due to Cartan \cite{cartan23}; for modern proofs, the reader is referred to Chapter 8 of \cite{sharpe97} or Sections 1.6.4--1.6.7 of \cite{cap-slovak-book09}. We note here that the normal conformal Cartan connection $\omega$ is always torsion-free, i.e., its curvature form satisfies $\Omega^{\omega} \in \Omega^2(\Gbdle;\palg)$.

\subsection{Reductive Cartan geometries}
\label{sec-background4} One class of Cartan geometries enjoying particularly nice properties are those of \emph{reductive type}: a type $(H,B)$ is said to be reductive if the Lie algebra $\halg$ has a decomposition
\begin{equation}\label{reductive decomposition}
\mathfrak{h} = \mathfrak{b} \oplus \mathfrak{n}
\end{equation}
into a direct sum of the Lie subalgebra $\mathfrak{b}$ and an $\mathrm{Ad}(B)$-invariant linear complement $\n$, i.e. such that $\mathrm{Ad}_H(B).\mathfrak{n} \subset \mathfrak{n}$. Usually, when we refer to a reductive Cartan geometry $(\pi: \mathcal{H} \rightarrow M,\eta)$ of type $(H,B)$, we will take a reductive decomposition (\ref{reductive decomposition}) to be fixed, although strictly speaking such a decomposition need not be unique. Note that the type $(G,P)$, as defined in Section \ref{sec-background3} for conformal geometries, is very far from being reductive (which is why we use different letters to denote the type). But Cartan geometries of reductive type will naturally occur when we study holonomy reductions corresponding to isotropy representations of certain symmetric spaces, so it is useful to collect some general properties of reductive Cartan geometries.

Central to the special properties enjoyed by reductive Cartan geometries is the decomposition $$\eta = \eta^{\mathfrak{b}} + \eta^{\mathfrak{n}}$$ of the Cartan connection $\eta \in \Omega^1(\mathcal{H};\mathfrak{h})$ given by projecting onto the factors in the decomposition (\ref{reductive decomposition}). In particular, from the $\mathrm{Ad}_H(B)$-invariance of that decomposition and equivariance of $\eta$, it follows that $\eta^{\mathfrak{b}} \in \Omega^1(\mathcal{H};\mathfrak{b})$ and $\eta^{\mathfrak{n}} \in \Omega^1(\mathcal{H};\mathfrak{n})$ are well-defined, $\mathrm{Ad}_H(B)$-equivariant one-forms.

In fact, it follows from the property of Cartan connections that $\eta^{\mathfrak{b}}$ is a $B$-principal connection, and $\eta^{\mathfrak{n}}$ is horizontal and determines a reduction of the frame bundle of $M$ to the structure group $\mathrm{Ad}_H(B) \subset \mathrm{GL}(\mathfrak{n})$. This reduction is given by using the maps $\psi_u: T_xM \rightarrow \mathfrak{h}/\mathfrak{b} \cong \mathfrak{n}$ as in (\ref{Cartan frames}) for a general Cartan geometry, which in the reductive case can be simplified to $\psi_u(X) = \eta^{\mathfrak{n}}(\widehat{X})$. The affine connection $\nabla^{\eta}: \Gamma(TM) \rightarrow \Gamma(T^*M \otimes TM)$ induced by the principal connection $\eta^{\mathfrak{b}}$ has torsion $T^{\eta} \in \Gamma(\Lambda^2T^*M \otimes TM)$ which is related to the torsion $\Theta^\eta$ of $\eta$ as a Cartan connection by the formula:
\begin{equation}\label{torsion0}
\psi_u(T^{\eta}(X,Y))=
\Theta^{\eta}(\widehat{X},\widehat{Y}) - \left[ \psi_u(X),\psi_u(Y) \right]_{\n},
\end{equation}
for any $\widehat{X}, \widehat{Y} \in T_u\mathcal{H}$ projecting to $X, Y \in TM$. For proofs of these facts, see, for example,   \cite[Appendix A]{sharpe97} or \cite[Section \ref{sec-orbit}]{lotta04}. We now prove a consequence of these properties which will be useful in studying the geometry induced by certain torsion-free reductive Cartan geometries:

\bs\label{reductivetorsion prop}
Let $(\pi: \mathcal{H} \rightarrow M,\eta)$ be a torsion-free Cartan geometry of reductive type $(H,B)$ and assume that the Lie algebra $\mathfrak{h}$ admits an $\mathrm{Ad}_H$-invariant non-degenerate metric $K: \mathfrak{h} \times \mathfrak{h} \rightarrow \R$ with respect to which the reductive decomposition $\mathfrak{h} = \mathfrak{b} \oplus \mathfrak{n}$ is orthogonal. Then the reductive Cartan geometry induces a canonical metric $g^{\eta}$ and metric affine connection $\nabla^{\eta}$ with torsion $T^{\eta}$ such that:

(i) The $(3,0)$-tensor given by contracting $T^{\eta}$ with $g^{\eta}$ is totally skew symmetric;

(ii) $\nabla^{\eta} T^{\eta} = 0$;

(iii) $\mathrm{Hol}(\nabla^{\eta}) \subset \mathrm{Ad}_H(B) \subset \mathrm{O}(\n,K)$.
\es

\bprf
As noted above, if we fix a point $x \in M$, then any choice of $u \in \mathcal{H}_x$ defines a linear isomorphism $\psi_u: T_xM \rightarrow \n$, and for any $b \in B$ we have the equivariance relation $$\psi_{u.b}(X) = \mathrm{Ad}(b^{-1}).\psi_u(X),$$ noting that this equivariance makes sense in the reductive setting because of $\mathrm{Ad}_H(B)$-invariance of the decomposition (\ref{reductive decomposition}). Now, from the orthogonality assumption $K(\mathfrak{b},\n) = 0$ it follows that $K_{\vert \n}$ defines a non-degenerate metric on $\n$, which is $\mathrm{Ad}_H(B)$-invariant since $K$ was assumed to be $\mathrm{Ad}_H$-invariant (in particular, $\mathrm{Ad}_H(B) \subset \mathrm{O}(\n,K)$). This allows us to define a metric $g^{\eta}$ on $M$ which is induced from the Cartan geometry, via $$g^{\eta}_x(X,Y) := K(\psi_u(X),\psi_u(Y)), \,\, \mathrm{for} \,\, \mathrm{any} \,\, u \in \mathcal{H}_x.$$

Since the maps $\psi_u: T_xM \rightarrow \n$ define a reduction of the frame bundle of $M$ to $\mathrm{Ad}_H(B) \subset \mathrm{GL}(\n)$, the $B$-principal bundle connection $\eta^{\mathfrak{b}} \in \Omega^1(\mathcal{H};\mathfrak{b})$ induces an associated affine connection $\nabla^{\eta}$ on $M$, with $\mathrm{Hol}(\nabla^{\eta}) \subset \mathrm{Ad}_H(B)$ by construction, so $\nabla^{\eta}$ is metric with respect to $g^{\eta}$ since $\mathrm{Ad}_H(B) \subset \mathrm{O}(\n,K)$.

It remains to verify properties (i) and (ii) for the torsion $T^{\eta}$ of this affine connection. For this, note that the expression (\ref{torsion0}) for the torsion $T^{\eta}$ shows, under the assumption that $(\mathcal{H},\eta)$ is torsion-free, that $T^{\eta} \in \Gamma(\Lambda^2T^*M \otimes TM)$ is induced by the alternating bilinear map $T^\n: \n \times \n \rightarrow \n$ given by $T^\n(u,v) := - [u,v]_{\n}$. So property (i) follows from the corresponding property of $T^\n$, which is computed as follows using $\mathrm{Ad}_H$-invariance of $K$ and orthogonality $K(\n,\mathfrak{b}) = 0$:
\begin{align*}
K(T^\n(u,v),w) &:= -K([u,v]_{\n},w) = -K([u,v],w) \\
&= K(v,[u,w]) = K(v,[u,w]_{\n}) \\
&=: -K(v,T^\n(u,w)).
\end{align*}

Finally, property (ii) follows from well-known properties of associated linear connections, from the fact that $\mathrm{Ad}_H(B) \subset \mathrm{Aut}(T^\n) \subset \GL(\n)$, while this last inclusion is straightforward to verify from the definition of $T^\n$.
\eprf
\subsection{Isotropy irreducible pseudo-Riemannian symmetric spaces}
\label{sec-background5}
Symmetric  spaces were classified by Cartan \cite{cartan26I} and Berger \cite{berger57} and their {\em pseudo-Riemannian} holonomy is equal to the isotropy group. First we collect some properties of irreducible pseudo-Riemannian symmetric spaces $(M,g)$ that will be used later.

\begin{satz} \label{Ricci nonflat} Let $(M,g)$ be an irreducible non-flat pseudo-Riemannian symmetric space, i.e., the isotropy group $\mathrm{Iso}(M)_x$ acts irreducibly on $T_xM$.
Then the Ricci tensor $\mathrm{Ric}_g$ is non-zero.
\end{satz}
Indeed, from a result of Nomizu \cite[16.1]{nomizu54} it follows that $\mathrm{Iso}(M)$ is either semisimple or $g$ is flat.
So $\mathrm{Iso}(M)$ must be semisimple. Then $M = G/H$ where the symmetric decomposition of $\mathfrak{g}:= Lie(G)$ is effective and minimal (see \cite{alekseevsky10}). Now the above proposition follows from \cite[Proposition 1]{alekseevsky10}.\\

The following proposition is a direct consequence of Propositions 1 and 2 of \cite{alekseevsky-cortes01}.

\begin{satz} \label{stabilizer} Let $(M,g)$ be an irreducible simply connected non-flat pseudo-Riemannian symmetric space and let $R$ be its curvature tensor.
Then the stabilizer $\mathrm{Aut}(R)$ of $R$ in the pseudo-orthogonal group $\O(T_xM,g_x)$ is the isotropy group $\mathrm{Iso}(M)_x$.
\end{satz}

\begin{bem} The above propositions were well-known to Cartan when he constructed his theory of Riemannian symmetric spaces.
For example, Proposition \ref{stabilizer} is interpreted as saying that an isometry $u$ of $(T_x M, g_x)$ can be extended to an isometry of $M$ fixing $x$ if and only if $u$ preserves the curvature tensor $R$.
\end{bem}

Now let $\g = \h \dsum \m$ be a symmetric decomposition, i.e., $\h \subset \g$ is a subalgebra, $[\m,\m] \subset \h$ and $[\h,\m] \subset \m$, with $\g, \h$ semisimple of non-compact type. Then the Killing form $K_{\g}$ restricts to a pseudo-Riemannian metric on $\m$, and restricting the adjoint representation of $\g$ to $\h$ defines the isotropy representation $\ad: \h \rightarrow \so(\m,K_{\g})$ at the Lie algebra level. In view of applications to conformal holonomy, we will be concerned with studying the subalgebras $$\b := \mathfrak{stab}_{\h}(\R S) = \{ X \in \h \mid \exists c_X\in \rr: [X,S]=c_XS \},$$ for $S \in \m$ a non-zero null vector with respect to $K_{\g}$.

 For the following result recall the notion of a Cartan involution $\theta$ of a semisimple Lie algebra $\h$: A Cartan involution of $\h$ is a Lie algebra automorphism $\theta$ which is involutive, $\theta^2 = \mathrm{Id}_{\h}$, and such that the restriction of the Killing form $K_{\h}$ to the $+1$-eigenspace of $\theta$ is negative-definite, and the restriction of $K_{\h}$ to the $-1$-eigenspace is positive-definite. In particular, since the eigenspaces are orthogonal, for any $0 \neq X \in \h$, we have $K_{\h}(X,\theta(X)) < 0$, so $\theta$-invariant subspaces are always non-degenerate with respect to $K_{\h}$.
We prove the following result, concerning distinguished cases where $\b$ enjoys nice properties:

\begin{satz} \label{reductive stabilizer}
Let $\g=\h \dsum \m$ be a symmetric decomposition with $\g$ and $\h$ semisimple of non-compact type, and let $S\in \m$ be a null vector and $\b:=\mf{stab}_\h(\R S)$. If $\h$ has a Cartan involution $\theta$ such that $\theta(\b)=\b$, then
\begin{enumerate}[(i)]
\item there is a reductive decomposition $\h=\b\+\n$ which is orthogonal with respect to the Killing form of $\g$, and  hence is naturally reductive, and
\item there is another null vector $\widehat{S} \in \m$ such that $K_{\g}(S,\widehat{S}) \neq 0$ and for $\widehat{\b} := \mf{stab}_{\h}(\R \widehat{S})$ we have
\begin{equation}\label{stab}
\widehat{\b} =  \b.
\end{equation}
\end{enumerate}
\end{satz}


\bbem
In  \cite[p. 207]{onishchik-vinberg3} a subalgebra $\b$ in a semisimple Lie algebra $\h$ that is invariant under a Cartan decomposition of $\h$ is called {\em canonically embedded}. In \cite[Theorem 3.6 in Chapter 6]{onishchik-vinberg3} it is shown that, for an algebraic subalgebra $\b$ of a real semisimple Lie algebra $\h$, this property is equivalent to $\b$ being reductive in the sense that $\ad(\b)$ is the tangent algebra of a reductive algebraic subgroup in $\GL(\h)$.
The proof uses the following fact (see
\cite[Theorem 2 in Chapter 4]{onishchik-vinberg90} or \cite[Chapter 1, Proposition 6.2]{onishchik-vinberg3}:
\bnum
\item[($\ast$)]
 Let $\mf{f}\subset \gl (V)$ be an algebraic linear Lie algebra over $\ccc$. Then $\mf{f}$ is the tangent algebra of a reductive algebraic linear group $F\subset \GL(V)$  if and only if the invariant scalar product $\trace (X\cdot Y)$ is non degenerate on $\mf{f}$.
 \enum
  This fact can also be used to establish the equivalence of our assumption that $\b$ is invariant under a Cartan involution with item (i) in Proposition  \ref{reductive stabilizer}. Indeed, applying ($\ast$) to $\mf{f}:=\b^\ccc$ and
  $V:=\g^\ccc$ then the trace form is given by the complexification $K_\g^\ccc$ of the Killing form $K_\g$ of $\g$ which is non-degenerate on $\h$ by Cartan's solvability criterion \cite[p. 68]{jacobson62}. Now, $K_\g$ is non-degenerate on $\b$ if and only if $K_\g^\ccc$ is non-degenerate on $\b^\ccc$. Observing that $\b$ is reductive if and only if $\b^\ccc$ is reductive and that in our situation $B$ and $H$ are defined as stabilizers and hence algebraic, yields the required equivalence. The notion of reductivity used in the statements  implies that the radical is equal to the center, but is stronger than that.
For our purposes it is sufficient to give a self contained proof of the weaker statement in Proposition \ref{reductive stabilizer} and to avoid subtleties in the notion of reductivity.
\ebem

\begin{proof}[Proof of Proposition \ref{reductive stabilizer}]
%
In order to prove the existence of another null line that is stabilised by $\b$ we apply
 the Karpelevich-Mostow Theorem (\cite{karpelevich53}, for the algebraic version we are using see \cite{mostow55}), which states that if $\h$ is a semisimple subalgebra of a semisimple Lie algebra $\g$ of non-compact type, then any Cartan involution $\theta$ of $\h$ extends to a Cartan involution $\widehat{\theta}$ of $\g$.
By the assumption, we have a Cartan involution $\theta$ of $\h$ such that $\theta(\b) = \b$. Let $\widehat{\theta}$ be a Cartan involution of $\g$ which extends $\theta$. Then $S$ decomposes as $S=S_++S_-$ into $(\pm 1)$-eigenvectors of $\widehat{\theta}$, and we define
%
 \[\hat S:=\hat\theta(S)=S_+-S_-\in \m.
 \]
Since $K_\g(S,S)=0$ we have that $S_\pm\not=0$ and hence $\hat S$ is linearly independent of $S$. Furthermore, since $S_+\bot S_-$, we have that
\[0=K(S,S)=K(S_+,S_+)+K(S_-,S_-)=K(\hat S,\hat S)\]
and $K(S,\hat S)=2K(S_+,S_+)<0$.

From the fact that $\widehat{\theta}$ is an involutive automorphism,
we get that $[X,S]=c_X S$ if and only if $[\theta (X), \hat S]=c_X \hat S$ for every $X\in \b$ and  a real constant $c_X$, which shows that
$\widehat{\b}:=\mf{stab}_\h(\hat S) = \theta(\b)$, which in turn equals $\b$ by $\theta$-invariance, showing \eqref{stab}.

In order to prove the first point, the invariance of $\b$ under $\theta$ implies the existence of one-forms $c$ and $\hat c$ of $\b$ such that
 $[X,S]=c(X) S$ and $[X,\hat S]=\hat c(X) \hat S$ satisfying
\[
c(X)\hat c(X) K(S,\hat S)=K([X,S],[X,\hat S]) =-K( S,[X,[X,\hat S]])=-\hat c(X)^2K(S,\hat S).\]
This  shows that $c=-\hat c$, since $K(S,\hat S)\not=0$, which in turn implies that $[X,S_\pm]=c(X) S_\mp$ for all $X\in \b$. When splitting $X\in \b$ as $X=X_++X_-$ with $X_\pm$ eigenvectors of $\theta$, we obtain that $[X_+,S_\pm]=0$ and $[X_-,S_\pm]=c(X) S_\mp $ and thus that $X_+\in \b$. This shows that for $X\in \b$ we also have $X_\pm\in \b$, and thus
\[ \b=
(\b\cap \h_+)\+^\bot(\b\cap \h_-)
\]
in which $\h_\pm$ denote the eigenspaces of $\theta $ in $\h$.
Since $K_\h$ and $K_\g$ are  negative-definite on $\h_+$ and positive-definite on $\h_-$ this shows that $K_\g$ is non-degenerate on $\b$. Now we get $\n:=\b^\bot$ as the reductive complement of $\b$ in $\h$, i.e., $\h=\b\+^\bot\n$ is naturally reductive.
\eprf
This proposition provides us with the main result of this section  which will be
 useful for proving Theorems \ref{sutheo} and \ref{sptheo} in Section \ref{applications}.

\btheo \label{Cartan invariant conformal holonomy} Let $\g = \h \dsum \m$ be a symmetric decomposition with $\g$ and $ \h$ semisimple of non-compact type, $K_\g$ the Killing form of $\g$,  and  let $G/H$ be the corresponding symmetric space. 
Let $(M,[g])$ be a conformal manifold of signature $(p,q)$ and 
suppose that \bnum[(i)]
\item $(M,[g])$ has a conformal holonomy reduction to the isotropy group $\Ad_G(H) \subset \SSO(\m,K_\g)\simeq\SSO(p+1,q+1)$,
\item \label{ii}
there is 
 is a null vector $S \in \m$  with stabilizer subgroup $B = \mathrm{Stab}_{H}(\R S)$ such that   the Lie algebra $\b$ of $B$ is invariant under a Cartan involution of $\h$. 
 \enum
Then the curved orbit $M_{S} \subset M$ corresponding to this $H$-orbit has a canonical metric $g_0$ and a canonical metric connection $\nabla^0$ with totally skew-symmetric, $\nabla^0$-parallel torsion $T^0$ and with holonomy contained in $\Ad_H(B) \subset \SSO(\h/\b)$. Moreover, if $(\n,K_{\g}\vert_{\n})$ and $(\widehat{\n},K_{\g}\vert_{\widehat{\n}})$ are homothetic, where $\n$ is the naturally reductive complement of $\b$ in $\h$ as in Proposition \ref{reductive stabilizer} and $\widehat{\n}$ is the $K_\g$-orthogonal complement to $\mathrm{span}(S,\widehat{S})$ in $\m$,  then the canonical metric $g_0$ is a representative of the conformal class $g|_{M_S}$.
\etheo

\bprf By Theorem \ref{cgh theorem}, there is a Cartan connection of type $(H,B)$ over $M_S$, which is a reduction of the canonical conformal Cartan connection of $(M,[g])$, and in particular is torsion-free. Because of  assumption \eqref{ii}, Proposition~\ref{reductive stabilizer} gives us  a decomposition $\h = \b \dsum \n$ which is naturally reductive with respect to the Killing form $K_\g$ of $\g$. Using 
 Proposition \ref{reductivetorsion prop}, this implies that
$M_S$ has a metric $g_0$ canonically induced  from $K_\g$ on $\h/\b\simeq\n$, and metric connection $\nabla^0$ with parallel, skew-symmetric torsion and $\mathrm{Hol}(\nabla^0) \subset \Ad_H(B)\subset \O(\n,K_\g|_\n)$, which proves the first statement.  Finally, recalling from 
Section \ref{sec-background3} that a representative in the conformal class is determined by 
pulling back 
$K_\g$ from $\g/\p\simeq \widehat\n$ to $T_xM$ via the isomorphism in \eqref{Cartan frames},
 the assumption that $(\n,K_{\g}\vert_{\n})$ and $(\widehat{\n},K_{\g}\vert_{\widehat{\n}})$ are homothetic
 shows that   the metric $g_0$  is in fact a representative of $[g_{\vert M_S}]$.
\eprf
Note that under the assumption that $\n$ and $\widehat\n$ are homothetic and hence have the same dimension, the orbit of an $S\in \m$ under  $H$ is open in the M\"obius sphere of $\m$.

\section{The orbit structure in the homogeneous models}
\label{sec-orbit}

\subsection{Semisimple symmetric spaces defined by ((hyper-)Hermitian-) scalar products}\label{secgroupdefs}
In the following, we will consider the isotropy representation of the semisimple, pseudo-Riemannian symmetric spaces $G/H$, where $G$  is given by
\[ G=\SSL_n\K,\ \ \text{ for $\K=\R$, $\C$, or the quaternions $\Q$}
\]
and the isotropy group $H$ is given by
\be
\SSO(p,q),&\text{if}& \K=\R,
\\
\SSU(p,q),&\text{if}& \K=\C,
\\
\SSp(p,q),&\text{if}& \K=\Q,
\ee
for $p+q=n$. For $\K=\R,\C$, $\SSL_n\K$ is the group of matrices with determinant one, while for $\K=\Q$ the special linear group $\SSL_n\Q$ is defined as the commutator group in $\GL_n\Q$. Given an identification of $\Q$ with $\C^2$ and the corresponding monomorphism of real algebras $\iota: \mathrm{Mat}_n\Q\hookrightarrow \Mat_{2n}\C$, $\SSL_n\Q$ is given as the preimage of matrices with determinant one (see for example \cite{aslaksen96} for a nice overview on quaternionic determinants).
The monomorphism $\iota$ can be given, for example, as
\be\iota\ :\  \mathrm{Mat}_n\Q&\hookrightarrow & \Mat_{2n}\C\\
U+V\cdot \vect{j}&\mapsto&\begin{pmatrix} U&-V\\\ol{V}&\ol{U}\end{pmatrix},
\ee
which satisfies $\iota(\ol{W}^\top)=\ol{\iota(W)}^\top$. We then have that
\[\iota(\SSL_n\Q)=
\SSU^*(2n):=\{A\in \SSL_{2n}\C \mid A\vect{J}_n=\vect{J}_n\ol{A}\},\]
where $\ol{A}$ denotes the complex-conjugated matrix and
\[\vect{J}_n=\begin{pmatrix}0&\1_n \\-\1_n&0\end{pmatrix}\in \GL_{2n}\R,\]
cf. \cite{helgason78}.
The Lie algebras of $\SSL_n\K$, for $\K=\R,\C$, are given as traceless matrices, whereas for $\K=\Q$ we have the real Lie algebra
\[\sl_n\Q:=\{ X+Y\, \vect{j}\mid X,Y\in \gl_n\C, \trace(X)+\trace(\ol{X})=0\}. \]
Now we define the isotropy group $H$ in $G$ as the invariance group of a ((hyper-)Hermitian-) scalar product
\[\la u,v\ra=\sum_{i=1}^p \ol{u^i}v^i- \sum_{j=p+1}^q \ol{u^j}v^j,\]
for $u,v\in \K^n$, which is anti-linear in the first slot. Here we consider $\Q^n$ as right vector space. For the standard basis in $\K^n$, $\la.,.\ra$ is given by the matrix
\[\1_{p,q}:=\begin{pmatrix}\1_p&0\\0&-\1_q\end{pmatrix}\in \GL_n\R.\]
For $\K=\R$ and $\K=\C$ we have
$H = \SSO(p,q)$ and $H = \SSU(p,q)$, whereas for $\K=\Q$ we have
\[
H\ =\ \SSp(p,q)
\ =\ \left\{A+B\,\vect{j}\mid \ol{A}^\top\1_{p,q}A+B^\top\1_{p,q}\ol{B}=\1_{p,q},\ B^\top\ipq  \ol{A}- \ol{A}^\top \1_{p,q}B=0\right\}.
\]
Note that $\iota(\SSp(p,q))\subset \SSU(2p,2q)$ when $\SSU(2p,2q)$ is written as
\[\SSU(2p,2q)=\left\{ A\in \sl_{2n}\C\mid \ol{A}^\top\mathbf{K}_{p,q}A=\mathbf{K}_{p,q}\right\},\]
where $\mathbf{K}_{p,q}=\begin{pmatrix}\1_{p,q}&0\\0&\1_{p,q}\end{pmatrix}$. With this realization of $\SSU(2p,2q)$ we have
\[
\iota(\SSp(p,q))= \SSU(2p,2q)\cap \SSp_n\C,\]
where the symplectic group $\SSp_n\C$ is defined as
\[
\SSp_n\C=\left\{A \in \GL_{2n}\C\mid A^\top \vect{J}_{p,q} A=\vect{J}_{p,q}\right\}
,\]
with the symplectic form $\vect{J}_{p,q}=\begin{pmatrix}0&\1_{p,q} \\ -\1_{p,q}&0\end{pmatrix}
$.

The Lie algebras $\h$ of the $H$'s are given as $\so (p,q)$ and $\su (p,q)$ and for $\K=\Q$  as
\[\sp(p,q):=
\left\{ X+\ipq Y\,\vect{j}\mid X\in \u(p,q),\, Y\in \gl_n\C \text{ symmetric}\right\}. \]
From the relation for the groups we get
\[\iota(\sp(p,q))=\su(2p,2q)\cap \sp_n\C,\]
where $\su(2p,2q)$ and $\sp_n\C$ again are defined with respect to $\vect{K}_{p,q}$ and $\vect{J}_{p,q}$.
This gives decompositions
 of $\g=\sl_n\K$ as a symmetric pair into
 $\g=\h\+\m^\h$ with
 \begin{equation}\label{m}
 \m^\h=\left\{X\in \sl_n\K\mid \ol{X}^\top\1_{p,q}=\1_{p,q}{X}\right\},\end{equation}
 for $\K=\R, \C$, i.e. for  $\h=\so(p,q)$ and $\h=\su(p,q)$, and with
 \begin{equation}\label{mhat}
 \m^{\sp(p,q)}=\{ X+\ipq Y\,\vect{j}\in \sl_n\Q \mid X\in \m^{\su(p,q)},\, Y\in \so_n\C \}\end{equation}
 in the quaternionic case. Furthermore, we have
 \blem \label{mlem}
 For $\mf{sym}(\vect{J}_{p,q}):=\left\{ W\in \gl_{2n}\C\mid W^\top \vect{J}_{p,q}=\vect{J}_{p,q}W\right\}$, then
 \[ \iota\left(\m^{\sp(p,q)}\right) =\m^{\su(2p,2q) }\cap \mf{sym}(\vect{J}_{p,q}).\]
 \elem
 \bprf The inclusion $\subset$ is verified  by a straightforward computation. For the other inclusion, we see that $W=\bmat X& -Y\\U&V\emat\in \m^{\su(2p,2q)}$ gives $X\in\m^{\su(p,q)}$ and $U= -\1_{p,q}\ol{Y}^\top \1_{p,q}$, whereas $W\in \mf{sym}(\vect{J}_{p,q})$ implies that $\ol{Y}^\top\1_{p,q}+\1_{p,q}Y=0$ and $X^\top\1_{p,q}=\1_{p,q}V^\top$. Together this implies that $W\in \iota(\m^{\sp(p,q)})$.
 \eprf
The decomposition $\g=\h\+\m^\h$ is invariant under the adjoint representation of $G$ when restricted to $H$. Therefore, the isotropy representation of $H$ is given by the adjoint representation of $G$ on $\m^\h$ restricted to $H$, i.e for $S\in \m^\h$ and $A\in H$, we have
 \begin{equation}\label{actionm} A( S):= \mathrm{Ad}_A(S)=A S A^{-1}.\end{equation}
The Killing form ${K}_\g$ of $\g$ is invariant under $\Ad_G$ and non-degenerate on $\m^\h$. Hence, we have that
 \[ \Ad_G(H)\subset \SSO(\m^\h \, , \,{K}_\g|_{\m^\h}).\]
 Note that, although in the case of the complex Lie algebra $\sl_n\C$ the Killing form $K_{\sl_n\C}$ is a complex bilinear form,  its restriction  to the real vector space $\m^{\su(p,q)}=\mathrm{i}\cdot \su(p,q)$ is real valued, in fact it is equal to $K_{\sl_n\C}(X,Y)=-K_{\su(p,q)}(\mathrm{i} X,\mathrm{i} Y)\in \rr$ for $X,Y\in \m^{su(p,q)}$.
 Furthermore,
we recall that for all the Lie algebras $\g$ and $\h$ in question, the Killing forms $K_\g$ and $K_\h$ are given as a real multiples of the trace form $(X,Y)\mapsto \trace(X\cdot Y)$, which allows us to determine the signature of $K_\g|_\m$. For convenience, we list the dimensions of $\h$, $\m^\h$ and the signature of the Killing form in Figure \ref{dimension table}.

\begin{figure}
\[\begin{array}{c|c|c|c}
\K&\dim\h&\dim \m^\h&\mathrm{sign }(K_\g |_{\m^\h})=(\sigma^+_{\K,p,q},\sigma^-_{\K,p,q})=(\text{no. of $+$'s, no. of $-$'s})  \\
\hline
&&&\\[-2mm]
\R&\einhalb (n-1)n& \einhalb (n-1)(n+2) &
\left( \tfrac{1}{2}\left(p(p+1)+q(q+1)-2\right),\ \  pq \right)\\
\C& n^2-1& n^2-1&\left( p^2+q^2-1,\, 2pq\right)\\
\Q& (2n+1)n&(2n+1)(n-1)& \left( p(2p-1)+q(2q-1)-1,\, 4pq\right)
\end{array}\]
\caption{Dimensions and signatures of the symmetric spaces $\SSL_n\K/H$}
\label{dimension table}
\end{figure}

Note also that for the Lie algebra $\h$ a Cartan involution is given by the transposition in $\sl_n\rr$ for $\h=\so(p,q)$ and by the conjugate transposition in $\sl_n\C$ and $\sl_{2n}\C$ for $\h=\su(p,q)$ and $\h=\sp(p,q)$, respectively.

In the remainder of Section \ref{sec-orbit}, we will analyze the orbit structure of the M\"obius sphere of $\m$ under the naturally induced $H$-action. Namely, for the null cone
\[ \mathcal{N} := \{ S \in \m \mid  K_\g(S,S) = 0\}\] in $\m$ and the projection \[\pi :\m \setminus \{0\} \to \mathbb{P}(\m)\]
onto the real projectivization of $\m$, the M\"obius sphere of $\m$ is defined as
\[
\mathbb{S}(\m):=
 \pi(\mathcal{N}).\]
 Since $H\subset \O(\m,K_\g|_\m)$ via the adjoint action as in \eqref{actionm}, in the same manner $H$ acts on $\cal N$ and on the M\"obius sphere $\mathbb{S}(\m)$,
\begin{equation}
\label{action}A([S]):=[ASA^{-1}], \text{ with } S\in \cal N, A\in H,\end{equation}
where we write $[S]=\pi(S)$ for brevity.
We define the stabiliser subgroup of $[S]$ as 
\begin{equation}\label{stabdef}\mathrm{Stab}_H([S]) := \{ A \in H : A([S])=[S] \}.\end{equation}
Note that, as $\dim \SS(\m)=\dim(\m)-2$, from Figure \ref{dimension table} we see that $\SSO(p,q)$ can only have open orbits in $\SS(\m^{\so(p,q)})$ if $n\le 3$, whereas for $\K=\C$ and $\K=\Q$ we always have $\dim(\h)>\dim (\SS(\m^\h))$. Note also that $n=2$ is not relevant for conformal holonomy, since for $\R$ and $\C$ the dimensions of $\m$ are too low, whereas $\sp(1,1)\simeq \so(4,1)$, which is the generic holonomy algebra of a conformal Riemannian $3$-manifold.

In the remainder of Section \ref{sec-orbit}, we will prove the following:

\btheo\label{orbit-theo}
Let $G/H = \SSL_n\K/H$ be one of the pseudo-Riemannian symmetric spaces defined above, and consider the natural $H$-action on the M\"obius sphere $\SS(\m)$ induced by the isotropy representation $\Ad_G: H \rightarrow \SSO(\m,K) \cong \SSO(\sigma^+_{\K,p,q},\sigma^-_{\K,p,q})$. Then the union of $H$-orbits of codimension $n-3$ is dense in $\SS(\m)$. In particular, for $\K = \R$ and all $n \geq 3$, the stabilizer subgroup 
$\mathrm{Stab}_H([S])$ 
is discrete for all $[S]$ in a dense subset $\SS(\m)_0 \subset \SS(\m)$. For $n=3$, the union of open $H$-orbits is a dense subset $\SS(\m)_0 \subset \SS(\m)$ for each case $\K=\R,\C$ or $\Q$, while for $n>3$ there are no open $H$-orbits. Finally, for any $[S] \in \SS(\m)_0$, there is a Cartan involution of $\h$ which leaves the stabilizer subalgebra $\b = \mf{stab}_{\h}([S])$ invariant.
\etheo

The proof proceeds by cases: first $\C$, then $\R$, then $\Q$. Note by Figure \ref{dimension table} that it suffices, in each of these respective cases, to prove that the  stabilizer subalgebra $$\mathfrak{stab}_{\h}(\R S) = \{ X \in \h : [X,S] = rS \,\, \mathrm{for} \,\, \mathrm{some} \,\, r \in \R \}$$ has real dimension $n-1$, $0$ and $3n$, respectively, for all null vectors $S$ in some dense subset $\mathcal{N}_0 \subset \mathcal{N}$. For the final claim of the Theorem, note that the map $\theta: X \mapsto -X^* := -\ol{X}^{\top}$, i.e. minus the conjugate-transpose, gives a Cartan involution of our Lie algebras $\h$ (indeed of $\g$) in all cases. We will see directly that the stabilizer subalgebras $\b$ are $\theta$-invariant, as part of the proofs determining their dimensions.

\subsection{Dense orbits: The special unitary case $\K=\C$}

\bs \label{unitary prop}
Let $n \geq 3$, and $G/H = \SSL_n \C / \SSU(p,q)$ as above, for $p,q\geq 1$. Let $\mathcal{N}_0 \subset \mathcal{N} \subset \m$ be the set of all null vectors $S \in \mathcal{N}$ which have mutually distinct eigenvalues $\lambda_1,\ldots,\lambda_n \in \C$. Then $\mathcal{N}_0$ is dense in $\mathcal{N}$ and, for all $S \in \mathcal{N}_0$ we have that the stabiliser in $\h$ of $S$ is conjugated to
$$\mathfrak{stab}_{\h}([S]) \simeq \{ \mathrm{diag}(z_1,\ldots,z_r,ix_1,\ldots,ix_{n-2r},-\ol{z_r},\ldots,-\ol{z_1}) \in \sl_n\C\},$$ for $z_1,\ldots,z_r \in \C$, $x_1,\ldots,x_{n-2r} \in \R$ and some $1 \le  r \leq n/2$,  with respect to a basis of eigenvectors for $S$. In particular, the stabilizer subalgebra has real dimension $n-1$ and it is invariant under the conjugate-transpose map $X \mapsto X^*$ in $\sl_n\C$. For $n=3$ we have $\mathrm{Stab}_H([S]) \cong \mathrm{U}(1) \times \mathrm{O}(1,1)$.
\es

\begin{proof} First we note that $\mathcal{N}_0$ is a dense open subset of $\mathcal{N}$. This follows, since $\mathcal{N}_0$ is the complement of the matrices in $\mathcal{N}$ whose characteristic polynomial has vanishing discriminant. The discriminant of the characteristic polynomial is polynomial in the entries of the matrix, and hence the complement of $\mathcal{N}_0$ is given as the zero set of an analytic function on $\mathcal{N}$, so it must either be all of $\mathcal{N}$ or have empty interior. And there certainly do exist null matrices $S \in \mathcal{N}_0$, as will be explained below.

To compute the explicit form of the Lie algebra $\mathfrak{b} := \mathfrak{stab}_{\h}([S])$, let $\{ u_1,\ldots,u_n \}$ be the basis of $\C^n$ consisting of eigenvectors of $S$. First note that, by the condition $K_{\g}(S,S) = \mathrm{tr}(S^2) = 0$, at least one of the eigenvalues, say $\lambda_1$, must be non-real. Furthermore, from the identities
\begin{align}
\lambda_i\la u_i,u_j \ra = \la Su_i,u_j\ra = \la u_i,Su_j\ra = \overline{\lambda_j}\la u_i,u_j\ra, \,\, 1 \leq i,j \leq n, \label{identities}
\end{align}
which follow from the defining equations of $\m \subset \sl_n\C$, it follows that $\la u_1,u_1\ra=0$; that, up to re-ordering and re-scaling we have $\la u_1,u_n\ra=1$ (by non-degeneracy of $\la.,.\ra$), $\lambda_n = \overline{\lambda_1}$, and $\la u_n,u_n\ra=0$; and that $\la u_1,u_j\ra = \la u_n,u_j\ra = 0$ for all $1 < j < n$. Clearly there is a maximal number $r \geq 1$ such that, after re-ordering, the eigenvalues $\lambda_1,\ldots,\lambda_r$ are all non-real and none of them are conjugate to each other. Using the identities (\ref{identities}) again, it follows for $1 \leq i \leq r$ that $\la u_i,u_i\ra = 0$; that there is a unique index $\nu(i)$, $r < \nu(i) \leq n$, such that $\la u_i,u_{\nu(i)}\ra \neq 0$; and that for this $\nu(i)$ we have $\lambda_{\nu(i)} = \overline{\lambda_i}$, and $\la u_{\nu(i)},u_j\ra=0$ for all $j \neq i$, $1 \leq j \leq n$. In particular, $r \leq \mathrm{min}(p,q)$. By re-ordering if necessary, we may take $\nu(i) = n-i+1$ for convenience, so the real eigenvalues are precisely $\lambda_{r+1},\ldots,\lambda_{n-r}$. Applying (\ref{identities}) to the corresponding eigenvectors shows, for all $r+1 \leq i \leq n-r$, that $\la u_i,u_j\ra = 0$ for all $j \neq i$. By non-degeneracy of $\la.,.\ra$, we must therefore have $\la u_i,u_i\ra \neq 0$ for all such $i$.

Summing up the above, we may assume after possibly re-ordering and re-scaling some of the eigenvectors $u_1,\ldots,u_n$, that the quadratic form of $\la.,.\ra$ has the following matrix form with respect to this basis:
\begin{equation}\label{a}
T_{p,q,r} := \begin{pmatrix}0 & 0 & R_r \\
0 & I_{p-r,q-r} & 0 \\
R_r & 0 & 0
\end{pmatrix},
\end{equation}
where $R_r = \mathrm{adiag}(1,\ldots,1)$ is the $r \times r$ matrix with $1$'s along the anti-diagonal and  $0$'s elsewhere. In general, by $\mathrm{adiag}(y_1,\ldots,y_r)$ be mean the anti-diagonal square $r \times r$ matrix having $y_1$ in the first row last column and $y_r$ in the last row first column, e.g., $\mathrm{adiag}(y_1,y_2) = \left(
                                                                                               \begin{array}{cc}
                                                                                                 0 & y_1 \\
                                                                                                 y_2 & 0 \\
                                                                                               \end{array}
                                                                                             \right)
$.

Now the form of $\mathfrak{b}$ can be calculated by simple linear algebra. Since $S$ is a diagonal matrix with mutually distinct diagonal entries, the equation $[X,S] = rS$ implies that all matrices $X \in \mathfrak{b}$ must be diagonal and hence $r=0$. Therefore,
\begin{align*}
\mathfrak{b} &= \{ X \in \mathfrak{su}(\C^n,\la.,.\ra) : X = \mathrm{diag}(X_1,\ldots,X_n) \}\\
&= \{ \mathrm{diag}(X_1,\ldots,X_n) : \ol{X}^{\top} T_{p,q,r} + T_{p,q,r} X = 0 \,\, \mathrm{and} \,\, \mathrm{tr}(X) = 0 \}.
\end{align*}
Using the form (\ref{a}) shows that the set of diagonal matrices $X \in \gl_n \C$ satisfying $\ol{X}^{\top} T_{p,q,r} + T_{p,q,r} X = 0$ have the form claimed in the proposition. In particular, this set has real dimension $n$, and therefore, since $\b \subset \sl_n\C$ and the trace of any matrix $X$ of that form is purely imaginary, $\mathrm{dim}(\mathfrak{b}) = n-1$.

Using a basis $\{u_1,\ldots,u_n\}$ in which $\la.,.\ra$ has the form (\ref{a}) also allows us to see that the set $\mathcal{N}_0$ is non-empty: For example, if we suppose $p \geq q \geq 1$, then we may take $$S = \mathrm{diag}(\mu_1,\ldots,\mu_q, \lambda_1,\ldots,\lambda_{p-q}, \overline{\mu_1},\ldots,\overline{\mu_q}),$$
with the $\mu_j = a_j + ib_j$ mutually distinct, non-real numbers, and the $\lambda_j$ mutually distinct, real numbers. Then $S \in \mathcal{N}_0$ is equivalent to the set of equations,
\begin{align}
2\sum_{j=1}^q a_j + \sum_{j=1}^{p-q} \lambda_j &= 0 \ \ (\mathrm{tr}(S) = 0); \label{trS=0}\\
2\sum_{j=1}^q (a_j^2-b_j^2) + \sum_{j=1}^{p-q} \lambda_j^2 &= 0 \ \ (\mathrm{tr}(S^2) = 0). \label{trS^2=0}
\end{align}
And the desired solutions to the above equations exist. This can be shown either by observing that the solution space of the simultaneous equations (\ref{trS=0}) and (\ref{trS^2=0}) forms a submanifold of positive dimension in the real parameters $\{a_1,\ldots,a_q,b_1,\ldots,b_q,\lambda_1,\ldots,\lambda_{p-q}\}$; or by an elementary direct construction of a solution.

The invariance under the conjugate-transpose map, $\b^* = \b$, follows immediately. Finally, for $n=3$, we have a complex basis $(u_1,u_2,u_3)$ of $\C^3$ such that $S$ is of the form $\mathrm{diag}(\mu,\lambda,\ol{\mu})$, with $\mu \notin \R$, $\lambda=-2\mathrm{Re}(\mu)$ and the scalar product $\la.,.\ra$ is of the form $R_3 = \mathrm{adiag}(1,1,1)$. For $A \in \mathrm{Stab}_H([S])$, the defining relation $ASA^{-1}=cS$ (for any $c \in \R^*$), or equivalently $AS=cSA$, therefore implies that $Au_1, Au_2, Au_3$ are eigenvectors of $S$ for the eigenvalues $\mu/c, \lambda/c$ and $\ol{\mu}/c$, respectively. But since $\mu, \lambda$ and $\ol{\mu}$ are distinct, this implies $c=1$ and $Au_i = \alpha_i u_i$ for $i=1,2,3$ and $\alpha_i \in \C$. Using the forms of $A$ and $\la.,.\ra$, one computes directly that $A\in\SSU(2,1)$ if and only if
\[
 A = \begin{pmatrix}r\e^{\frac{\mathrm{i}\vf}{2}}&0&0 \\ 0& e^{-\mathrm{i}\vf}&0
 \\
0&0&\frac{1}{r}\e^{\frac{\mathrm{i}\vf}{2}}
\end{pmatrix}
\]
for some $\vf, r \in \R, r \neq 0$, giving the isomorphism to $\U(1)\times \mathrm{O}(1,1)$ as claimed.
\end{proof}

\subsection{Dense orbits: The orthogonal case $\K=\R$}\nopagebreak[4]
\bs \label{real prop}
Let $n \geq 3$, and $G/H = \SSL_n \R / \SSO(p,q)$ for $p,q\geq 1$. Let $\mathcal{N}_0 \subset \mathcal{N} \subset \m$ be the set of all null vectors $S \in \mathcal{N}$ which, considered as endomorphisms acting on $\C^n$, have mutually distinct eigenvalues $\lambda_1,\ldots,\lambda_n \in \C$. Then $\mathcal{N}_0$ is dense in $\mathcal{N}$ and, for all $S \in \mathcal{N}_0$ we have
$$\mathfrak{stab}_{\h}([S]) = \{0\}.$$
\es

\bprf
Let $S\in \cal N\subset \m^{\so(p,q)}$ be a real matrix acting on $\rrn=\mathrm{span}(e_1, \ldots , e_n)$. When we consider $S$ as acting on $\C^n=\rrn\+\mathrm{i}\rrn$ by complex linear extension, and $\C^n$ as equipped with the Hermitian form coming from the real scalar product on $\rrn$, we have that $S\in \m^{\su(p,q)}$.
Now
 \[\cal N_0=\{S\in \cal N \subset \m^{\so(p,q)}\mid S \text{ has $n$ pairwise distinct eigenvalues over $\C$}\}
 \]
and we can use the results of the previous section.
We fix a complex basis $$\{v_1, \ldots , v_r, u_1, \ldots , u_{n-2r}, w_r, \ldots , w_1\}$$ of $\C^n$ in which $S$ is of the form
\[S=\mathrm{diag}(\mu_1, \ldots , \mu_r, \lambda_1, \ldots , \lambda_{n-2r}, \ol{\mu}_r, \ldots , \ol{\mu}_1),
\]
as given by the proof of Proposition \ref{unitary prop}.
Now denote by $\ol{ v}$ the conjugation on $\C^n$ induced by $\rrn\subset \C^n$.
 Clearly we have that $\ol{Sv}=S\ol{v}$ for each $v\in \C^n$. Indeed, for $v=\sumi (a^i+\mathrm{i}b^i)e_i\in \C^n$ with $a^i,b^i\in \rr$ we get
 \[
 \ol{S(v)}= \sumi \left(a^i S(e_i)-\mathrm{i} b^iS(e_i)\right)=S(\ol{v}),\]
 since $S(e_i)\in \rrn$. Applying this to the eigenbasis
  gives
  \be
  S(\ol{v}_i)&=&\ol{\mu}_i\ol{v}_i
  \\
  S(\ol{u}_k)&=&\lambda_k\ol{u}_k.
  \ee
Since the eigenvalues of $S$ are pairwise distinct,  this shows
   that $w_i=\ol{v}_i$ and $\ol{u}_k=u_k$. Hence, the vectors
 \be
 \{x_i := v_i+\ol{v}_i, y_j := \mathrm{i}(v_j-\ol{v}_j), u_k\},
 \ee
for $i,j=1, \ldots r$ and $k=1,\ldots , n-2r$, form a real basis of $\rr^n\subset \C^n$ in which the scalar product is diagonal with $\pm 2$ on the diagonal.

When we consider $X\in \so(p,q)$ as acting on $\C^n$, we get that $X\in \su(p,q)$. We have seen that the relation $[X,S]=cS$ for $c\in \rr$ implies that, in the eigenbasis for $S$, $X$ is given as
\[X=\mathrm{diag}(z_1, \ldots , z_r, \mathrm{i}s_1, \ldots , \mathrm{i}s_{n-2r}, -\ol{z}_r, \ldots , -\ol{z}_1)
\]
for $z_i\in \C$ and $s_k\in \rr$. But, from the deduced form of $X$, we get that
\be
 X(x_i)&=& \mathrm{i}( b_ix_i-a_i y_i )
 \\
 X(u_k)&=&\mathrm{i}s_ku_k,
 \ee
with $z_i=a_i+\mathrm{i}b_i$ the complex eigenvalues of $X$. This
is a contradiction to the invariance of $\rrn$ under $X$ unless $z_i=s_k=0$.
   Hence, the stabiliser of $S\in \cal N_0$ in $\so(p,q)$ is trivial.

 For the proof it remains to show that
 \[\cal N_0:=\{S\in \cal N \subset \m^{\so(p,q)}\mid S \text{ has $n$ pairwise distinct eigenvalues over $\C$}\}
 \]
 is dense in $\cal N\subset \m^{\so(p,q)}$. This follows as in the proof of Proposition \ref{unitary prop}, noting again that $\cal N_0$ is non-empty, since every matrix which has $2\times 2$ matrices of the form
\[A_i:= \begin{pmatrix} a_i&b_i\\-b_i&a_i\end{pmatrix},\ \text{ with }a_i\in \rr,b_i\in \rr^*,\]
and real numbers $c_1, \ldots c_{n-2r}$ on the diagonal, with respect to the above basis, is in $\cal N_0$ as long as the $A_i$'s and the $c_i$'s are mutually distinct, and the $a_i, b_i, c_i$ satisfy equations analogous to (\ref{trS=0}) and (\ref{trS^2=0}). \eprf

\begin{bem}
With similar computations as in the proofs of Propositions \ref{unitary prop} and \ref{real prop} it is possible to show that in both cases the open orbit is unique for $n=3$ and $\K = \R$ or $\C$. It follows that if $\det(S'), \det(S) \neq 0$ then $H.[S] = H.[S']$, for $[S], [S'] \in \SS(\m^{\so(2,1)})$, respectively for $[S], [S'] \in \SS(\m^{\su(2,1)})$. Moreover, in the case $\K = \R$ we were able to find an explicit description of all the $H$-orbits: the unique open orbit consists of the image in $\SS(\m^{\so(2,1)})$ of the invertible null matrices in $\m^{\so(2,1)}$; there is one orbit of codimension one, given by the image of all two-step nilpotent matrices in $\m^{\so(2,1)}$; and there is one orbit of codimension two, given by the image of one-step niplotent matrices in $\m^{\so(2,1)}$. We did not take the time to attempt the corresponding computations to find explicit descriptions of the $H$-orbits for $\K = \C$ or $\Q$, because these descriptions were not needed for the main applications in the paper.
\end{bem}

\subsection{Dense orbits: The symplectic case $\K = \Q$.}\label{sympsec}
Recall that we identify \[\sl_n\Q\ \simeq \ \su^*(2n)\ =\ \left\{\begin{pmatrix} X & -Y \\ \ol Y&\ol X\end{pmatrix}\mid X, Y\in \gl_n\C, \mathrm{tr}(X)+\mathrm{tr}(\ol X)=0\right\}
\]
and under this identification
\[
\sp(p,q)=
\left\{\begin{pmatrix} X & -Y \\ \ol Y&\ol X\end{pmatrix}\mid X\in \u(p,q), Y\in \gl_n\C: Y^\top\1_{p,q}-\1_{p,q}Y=0\right\}.
\]
Then $\m:=\m^{\sp(p,q)}$ is given as
\[\m:=
\left\{\begin{pmatrix} X & -Y \\ \ol Y&\ol X\end{pmatrix}\mid X\in\m^{\su(p,q)}, Y\in \gl_n\C: Y^\top\1_{p,q}+\1_{p,q}Y=0\right\}.
 \]
Furthermore, we have
\begin{equation}\label{msp}\m =\m^{\su(2p,2q) }\cap \mathfrak{sym}(\vect{J}_{p,q}),\end{equation}
for $\mathfrak{sym}(\vect{J}_{p,q})$, $\vect{J}_{p,q}$ and $\su(2p,2q)$ as described in Section \ref{sec-orbit}.1, cf. Lemma \ref{mlem}. We write $\la .,.\ra$ for the metric, given by $\mathbf{K}_{p,q}$, which determines $\su(2p,2q)$ and $\m^{\su(2p,2q)}$ as subspaces of $\sl_{2n}\C$, and $\omega$ for the symplectic form, given by $\vect{J}_{p,q}$, which determines $\sp_n\C$ and $\mathfrak{sym}(\vect{J}_{p,q})$ as subspaces of $\sl_{2n}\C$.

Now we consider the Jordan canonical form of elements of $\m$. Recall the following result, due to Wiegmann (\cite{wiegmann55}, see \cite{zhang97} for an overview and Corollary 6.3 therein): The Jordan canonical form of a complex matrix
\[Z:=\begin{pmatrix} X & -Y \\ \ol Y&\ol X\end{pmatrix}\]
with $X$ and $Y$ being complex $n\times n$ matrices is given by
\[\hat J=\begin{pmatrix} J & 0 \\ 0 &\ol J\end{pmatrix},\]
where $J$ is the Jordan normal form of some complex matrix. In particular, all the Jordan blocks come in pairs with complex conjugate eigenvalues. We call the eigenvalues of $J$ the generalized eigenvalues of $Z$. Furthermore, the matrix $B\in \GL_{2n}\C$ such that $B^{-1}ZB=\hat J$ is an element in $\iota(\GL_n\Q)$, i.e. of the form
\[B=\begin{pmatrix} P & -Q \\ \ol{Q} &\ol P\end{pmatrix}.\] In this section we will prove:

\bs \label{symplectic prop}
Let $n \geq 3$, and $G/H = \SSL_n \Q / \SSp(p,q)$ for $p,q\geq 1$. Let $\mathcal{N}_0 \subset \mathcal{N} \subset \m^{\sp(p,q)} \subset \m^{\su(2p,2q)}$ be the set of all null vectors $S \in \mathcal{N}$ which have mutually distinct generalized eigenvalues $\lambda_1,\ldots,\lambda_n \in \C$. Then $\mathcal{N}_0$ is dense in $\mathcal{N}$ and, for all $S \in \mathcal{N}_0$ we have a basis in $\C^{2n}$ in which the stabilizer subalgebra $\b = \mathfrak{stab}_{\h}([S])$ has the form
\[
\left\{ \begin{pmatrix}X&-Y\\\ol{Y}&\ol{X}\end{pmatrix}\left|
\begin{array}{l}
X=\mathrm{diag} (z_1, \ldots z_r, \mathrm{i}x_1,\ldots , \mathrm{i} x_{n-2r}, -\ol{z}_r, \ldots , -\ol{z}_1)\\
Y=
\begin{pmatrix}
0&0&\mathrm{adiag} (y_1, \ldots y_r)\\
0&\mathrm{diag}(y_{r+1}, \ldots , y_{n-r})&0
\\
\mathrm{adiag} (y_{n+1-r}, \ldots y_n)&0&0\end{pmatrix}
\\
y_i\in \C, i=1, \ldots, n, z_j\in \C, j=1,\ldots, r, x_k\in\rr, k=1, \ldots , n-2r
\end{array}
\right.
\right\},
\]
where $1\le r\le \frac{n}{2}$  and {\em adiag} denotes the anti-diagonal matrix, e.g., $\mathrm{adiag}(y_1,y_2) = \left(
                                                                                               \begin{array}{cc}
                                                                                                 0 & y_1 \\
                                                                                                 y_2 & 0 \\
                                                                                               \end{array}
                                                                                             \right)$.
 In particular, $\b$ is isomorphic as real Lie algebra to
 \[
 \b= \underbrace{\sp(1)\+\ldots \+ \sp(1)}_{n-2r\text{ times}} \+ \underbrace{\sl_2\C\+\ldots \+\sl_2\C}_{r\text{ times}},
 \]
has real dimension $3n$, and is invariant under the conjugate-transpose map $Z \mapsto Z^*$ of $\sl_n\Q \simeq \su^*(2n) \subset \sl_{2n}\C$, i.e. $\b^* = \b$. For $n=3$, the stabilizer subgroup $B = \mathrm{Stab}_H([S])$ is isomorphic to $\SSp(1) \times \SSL_2\C$.
\es

\bprf
Let us fix some notation. Let $T:=T_{n-r,r,r}$ denote the $n \times n$ matrix with the form of (\ref{a}). Furthermore, denote by $Q$ a matrix with (arbitrary) non-zero complex entries in precisely the positions where the matrix $T$ has $\pm 1$. Then we have:

\blem
Let $S\in \cal N_0$. Then there is a basis of $\C^{2n}$ of eigenvectors of $S$ such that the scalar product $\la.,.\ra$ and the symplectic form  $\omega$ are given, respectively, by
\[
\begin{pmatrix}
Q&0\\
0&Q
\end{pmatrix}\text{ and }
\begin{pmatrix}
0&T\\
-T&0
\end{pmatrix}
\]
in this basis.
\elem

\bprf
Let $(v_1, \ldots , v_n, w_1, \ldots w_n)$ be a basis of eigenvectors of $S \in \mathcal{N}_0$. From \eqref{msp} we know that
\[S\in \m=
\{A\in \gl_{2n}\C\mid \la Ax,y\ra=\la x,Ay\ra \text{ and }\omega(Ax,y)=\omega(x,Ay)\}\]
and thus
\begin{eqnarray}
(\lambda_i-\ol{\lambda_j})\la v_i,v_j\ra\ =\ (\lambda_i-\ol{\lambda_j})\la w_i,w_j\ra &=&0\label{one}
\\
(\lambda_i-{\lambda_j})\la v_i,w_j\ra&=&0\label{two}
\\
(\lambda_i-{\lambda_j})\omega(v_i,v_j)\ =\ (\lambda_i-{\lambda_j})\omega(w_i,w_j) &=&0\label{three}
\\
(\lambda_i-\ol{\lambda_j} )\omega(v_i,w_j)&=&0\label{four}
\end{eqnarray}
for $1\le i,j\le n$. Since $\lambda_i\not= \lambda_j$, \eqref{three} implies that
  \[  \omega(v_i,v_j)\ =\ \omega(w_i,w_j) \ = \ 0 \text{ for $ 1\le i,j\le n$,}\]
and (\ref{two}) implies that
  \[ \la v_i,w_j\ra=0  \text{ for $1\le i\not= j\le n$.}\]
Again, since $S^2$ has no trace, one of the $\lambda_i$'s must be non-real, let's say $\lambda_1$. Thus \eqref{four} implies that $\omega(v_1,w_1)=0$. Hence, as $\omega$ is non-degenerate, we can assume that $\omega(v_1,w_n)=1$. This implies that
  \[\lambda_1=\ol\lambda_n.\]
As the $\lambda_i$'s are pairwise distinct we get
\be
  \lambda_1&\not=& \ol\lambda_j,\ \text{ for }j=1, \ldots, n-1\\
  \lambda_n&\not=&\ol\lambda_j,\ \text{ for }j=2, \ldots, n
  \ee and therefore
  \barr{rcccl}
  \la v_1, v_j\ra&=&\omega(v_1,w_{j}) &=&0,\ \text{ for }j=1, \ldots, n-1\\
\la v_n,v_j\ra&=&\omega(v_n,w_j) &=&0,\ \text{ for }j=2, \ldots, n.\\
\earr
Hence, if we have $\lambda_1, \ldots , \lambda_r\not\in \rr$, in a similar way we get
\be
\lambda_{n+1-i}&=&\ol\lambda_i,\    \text{ for }i=1, \ldots ,r
\ee
and
\be
\la v_i,v_{k}\ra\
=\
\la w_i,w_{k}\ra\
=\ \omega(v_i,w_{k})&=&0,\ee
 for $ i\in \{1, \ldots r\}\cup\{ n-r+1, \ldots n\}, k\not= n+1-i$. Furthermore, for $k=r+1, \ldots , n-r$ the $\lambda_k$'s are real and we have
 \[
 \la v_k,v_l\ra = \omega(v_k,w_l) =0
 \]
for $r+1\le k\not= l\le n-r$. Hence, the symplectic form is represented by the desired matrix
$\begin{pmatrix}
0&T\\
-T&0
\end{pmatrix}
$.
Furthermore, the only possible non-vanishing terms for $\la.,.\ra$ in this basis are
 \barr{rl}
 \la v_i,w_i\ra, & \text{ for }1\le i\le n
 \\
 \la v_i,v_{n+1-i}\ra,\  \la w_i,w_{n+1-i}\ra, & \text{ for } i\in \{1, \ldots r\}\cup\{ n-r+1, \ldots n\}, k\not= n+1-i
\\
 \la v_k,v_k\ra,\  \la w_k,w_{k}\ra, & \text{ for } r+1\le k\le  n-r
 \earr
In order to change the basis to achieve the required form for $\la.,.\ra$ we note that
\[ S=\mathrm{diag}( \lambda_1, \ldots ,\lambda_r,\lambda_{r+1},\ldots, \lambda_{n-r}, \ol\lambda_{r},\ldots ,\ol\lambda_1, \ol\lambda_1,\ldots \ol\lambda_{r}, \lambda_{r+1}, \ldots,\lambda_{n-r},
\lambda_r, \dots, \lambda_1),
\]
which shows that $S$ has $n$ two-dimensional eigenspaces, $n-2r$ many for the real eigenvalues $\lambda_{r+1},\ldots, \lambda_{n-r}$,
\[U_k:=\mathrm{span}(v_k,w_k),\ k=r+1, \ldots, n-r,\]
and $2r$ many for  the complex eigen values
$\lambda_1,\ldots, \lambda_r$ and $\ol{\lambda_1},\ldots , \ol\lambda_r$,
\[ V_i=\mathrm{span}(v_i,w_{n+1-i})\text{ and } W_i=\mathrm{span}(v_{n+1-i},w_{i})
\]
for $i=1, \ldots, r$. Note that $U_k$ is orthogonal to $U_l$ for $k\not=l$ and orthogonal to $V_i$ and $W_i$.
This allows us to change the basis within the $U_k$'s in a way that $\la.,.\ra$ is diagonal on $U_k$.
Note that the diagonal does not have to be $\pm 1$. We can only use base change matrices
with determinant one in order to preserve the standard symplectic form.
Furthermore, $V_i$ and $W_i$ are totally isotropic, and $V_i\perp V_j$ and $W_i\perp W_j$ for $i\not=j$. The symplectic form  $\omega$ on $V_i$ and $W_i$ is given as the standard one $\vect{J}_2$. On $V_i\+W_i$, the scalar product and the symplectic form are given as
\[
\begin{pmatrix}0&H\\  H^\top&0\end{pmatrix},\ \
\begin{pmatrix}\vect{J}_2&0\\ 0&\vect{J_2}\end{pmatrix},
\]
where $H$ is an invertible $2\times 2$-matrix. Changing the basis of $V_i$ and $W_i$ by means of invertible matrices $A$ and $B$, respectively, yields the following matrices for
the scalar product and the symplectic form
\[
\begin{pmatrix}0&A^\top HB\\ B^\top H^\top A&0\end{pmatrix},\ \
\begin{pmatrix}\det(A)\vect{J}_2&0\\ 0&\det(B)\vect{J_2}\end{pmatrix}.
\]
This allows us to diagonalize $H$ and get the desired form for $\la.,.\ra$.
Namely, we can guarantee that the matrix of $\la.,.\ra$ has the required form $\mathrm{diag}(Q,Q)$ as in the Lemma by taking $A,B \in \mathrm{SL}_2\C$ such that $A^\top HB$ is a constant multiple of $\mathbf{1}_2$, e.g. take $A=\mathbf{1}_2$ and $B = \sqrt{\mathrm{det}(H)}H^{-1}$.
\eprf

Now we can prove Proposition \ref{symplectic prop}. We set ${\lambda}:=\mathrm{diag}(\lambda_1,\ldots,  \lambda_n)$ and fix $S=\begin{pmatrix}\lambda&0\\
0&\ol\lambda\end{pmatrix}\in \cal N_0$.
Then, for $U=\begin{pmatrix}X&-Y\\ \ol{Y}&\ol{X}\end{pmatrix} \in \sp (p,q)$ the relation
$[U,S]=c S$ for a real number $c$ amounts to the relations
\be
\lambda Y-Y\ol{\lambda}&=&0
\\
\lambda X-X\lambda&=&c\lambda.
\ee
These relations on one hand imply that $X=\mathrm{diag}(z_1, \ldots, z_n)$, with $z_i\in \C$ and on the other, due to the conditions on the $\lambda_i$'s that
\[
Y=
\begin{pmatrix}
0&0&\mathrm{adiag} (y_1, \ldots y_r)\\
0&\mathrm{diag}(y_{r+1}, \ldots , y_{n-r})&0
\\
\mathrm{adiag} (y_{n+1-r}, \ldots, y_n)&0&0\end{pmatrix}
\] with
$y_i\in \C$.

Then a straightforward computation shows that the invariance of the symplectic form $\omega$ under such a matrix
$U=\begin{pmatrix}X&-Y\\\-\ol{Y}&\ol X\end{pmatrix}\in \mathfrak{stab}(S)$, i.e.,
\[
U^\top\begin{pmatrix}
0&T\\
-T&0
\end{pmatrix}
+ \begin{pmatrix}
0&T\\
-T&0
\end{pmatrix}U=0,
\]
poses no further conditions on $Y$, but forces
\be
z_k+\ol{z_k}&=&0,\ \text{ for }k=r+1, \ldots, n-2r
\\
z_i+\ol{z}_{n+1-i}&=&0,\ \text{ for }i=1, \ldots, r.
\ee
Then one computes easily that  the scalar product $\la.,.\ra$ is invariant under such a matrix $U$, i.e., that
\[
\ol{U}^\top
\begin{pmatrix}
Q&0\\
0&Q
\end{pmatrix}
+ \begin{pmatrix}
Q&0\\
0&Q
\end{pmatrix}U=0
\]
holds true.

Now we verify the Lie algebra structure of $\b$:
We  set $\ell:=n-2r$, and represent an element of $\b$ by a pair $(X,Y)$  with $X$ and $Y$ of the form given in the statement of the Proposition. When writing the elements of $\b$ in this as pairs of matrices and denoting by $\e_i$ a standard basis vector in $\rr^n, \rr^r$ or $\rr^\ell$, we claim that the $\sp(1)$ summands are given as
\[
\b_k:=\left\{\left.\left(X_k:= \mathrm{diag}(\mathrm{i} a\e_{r+k} ),
Y_k:= \begin{pmatrix}0&0&0\\0&\mathrm{diag} (z \e_k)&0
\\0&0&0 \end{pmatrix} \right)\right| a\in \rr, z\in \C\right\}
 \]
 for $k=1, \ldots ,\ell$, and that the $\sl_2\C$ summands are given as
 \[
\mathfrak{s}_i:=\left\{\left. \left( X_i:= \mathrm{diag}(z\e_{i} - \ol{z}\e_{n+1-i}  ),Y_i:=
 \begin{pmatrix}0&0&\mathrm{adiag} (x\e_i)\\
 0&0&0
\\ \mathrm{adiag} (y\e_{2r+1-i}) &0&0 \end{pmatrix} \right)\right| x,y,z\in \C\right\}
 \]
for $i=1, \ldots , r$. Clearly, these spaces commute with each other and it  is a straightforward computation to check that they   enjoy the commutation relations of $\sp(1)$ and $\sl_2\C$. 

Finally, for the group isomorphism $B \simeq \SSp(1) \times \SSL_2\C$  for $n=3$ note 
that $ \SSp(1) \times \SSL_2\C$ can be embedded into $\SSp(2,1)$. Explicitely,
assigning to a unit quaternion $(u+v\vect{j})\in \SSp(1)$ matrices 
\[U:=\begin{pmatrix}1&0&0\\0&u&0\\0&0&1\end{pmatrix},\  V:=\begin{pmatrix}0&0&1\\0&v&0\\ 1&0&0\end{pmatrix},\]
a straightforward computation shows that $\begin{pmatrix}U&-V\\ \ol{V}&\ol{U}\end{pmatrix}$ is in $H\simeq \SSp(2,1)$ (see also  Appendix \ref{appA} for explicit formulae). Similarly, to any matrix $\begin{pmatrix}a&b\\ c&d\end{pmatrix}\in \SSL_2\C$ we can assign matrices 
\[A:=\begin{pmatrix}\overline{a}&0&0\\0&1&0\\0&0&d\end{pmatrix},\  C:=\begin{pmatrix}0&0&b\\0&1&0\\ -c&0&0\end{pmatrix},\]
such that 
$\begin{pmatrix}A&-C\\ \ol{C}&\ol{A}\end{pmatrix}$ is in $H\simeq \SSp(2,1)$. This gives an embedding of $\SSp(1) \times \SSL_2\C$ into $H\simeq \SSp(2,1)$, for which the Lie algebra of the image is given as $\b$. Hence, the stabiliser $B$ is isomorphic to $\SSp(1) \times \SSL_2\C$.
\eprf

\section{Applications to conformal holonomy}
\label{applications}

In this section, we apply the facts from Section \ref{sec-orbit} about orbits in the homogeneous models, and Theorem \ref{cgh theorem} of \cite{cgh11}, to study the geometry induced by conformal holonomy reductions to the isotropy subgroups
\begin{align*}
\Ad_{\SSL_3\R}(\SSO(2,1)) &\subset \SSO(3,2);\\
\Ad_{\SSL_3\C}(\SSU(2,1)) &\subset \SSO(4,4);\\
\Ad_{\SSL_3\Q}(\SSp(2,1)) &\subset \SSO(6,8).
\end{align*}
In particular, we will prove Theorems \ref{gensotheo}, \ref{sutheo} and \ref{sptheo}, showing ``essentially'' (i.e. after possibly restricting to an open dense subset of the conformal manifold) that these isotropy representations are not geometrically realizable as conformal holonomy groups.

\subsection{Proof of Theorem \ref{gensotheo}} Let $(M,[g])$ be a conformal manifold of signature $(2,1)$
and let $H \subset \SSO(3,2)$ denote the isotropy subgroup $H=\Ad_{\SSL_3\R}(\SSO(2,1))$. As in Section \ref{sec-background3}, we denote $G:=\SSO(3,2)$ and let $P \subset G$ be the parabolic subgroup stabilizing some null ray in $\R^{3,2}$. By Proposition \ref{real prop}, the union of the  induced $H$-orbits $H(gP)$ having (maximal) dimension $3 = \mathrm{dim}_{\R}(H)=\dim \SS^{2,1}$, are dense in the homogeneous model $G/P$ which is a double covering of $\SS^{2,1}$. Now, if $\mathrm{Hol}(M,[g]) \subseteq H \subset \SSO(3,2)$, then Proposition \ref{stabilizer} supplies us with a natural choice of parallel tractor $\Upsilon \in \Gamma(\bigotimes^4 \mathcal{T}^*)$ giving a reduction to $H \subset \SSO(3,2)$---namely, $\Upsilon$ corresponds to the pseudo-Riemannian curvature tensor $R \in \bigotimes^4(\R^{3,2})^*$ of the symmetric space $\SSL_3\R/\SSO(2,1)$. This we use to conveniently apply the results reviewed in Section \ref{sec-background2}, letting $\mathcal{O} := \rho(G).R$ denote the $G$-type of $\Upsilon$. First define the subset $M_0 \subset M$ to be the union of all curved orbits $M_{\ol{\alpha}} \subset M$ as in (\ref{curved orbit decomposition}), such that the $P$-type $\ol{\alpha} \in P \backslash \mathcal{O}$ corresponds to an $H$-orbit of dimension $3$ in $G/P$. It follows from Theorem \ref{cgh theorem}  that $M_0$ is dense in $M$. Furthermore, it follows that each point $x \in M_0$ lies in a curved orbit $M_{\ol{\alpha}}$ of dimension $3$ (by construction), and that the curved orbit $M_{\ol{\alpha}}$ carries a canonical Cartan geometry $(\mathcal{H}_{\ol{\alpha}} \rightarrow M_{\ol{\alpha}},\eta_{\ol{\alpha}})$ which is a reduction of the canonical conformal Cartan geometry $(\mathcal{G} \rightarrow M,\omega^{nc})$ to a Cartan geometry of type $(H,P_{\ol{\alpha}})$,
\be
\mathcal{H}_{\ol{\alpha}}& \stackrel{\iota}{\hookrightarrow}&\cal G
\\
\downarrow&&\downarrow
\\
 M_{\ol{\alpha}}&\subset &M
 \ee
 and that $\iota^*\omega^{nc} =\eta_{\ol{\alpha}}$.
 In particular, since $\omega^{nc}$ is torsion-free, so is the Cartan connection $\eta_{\ol{\alpha}}$, but this means that $\eta_{\ol{\alpha}}$ is flat, since the subgroup $P_{\ol{\alpha}}$ is discrete, as a consequence of Theorem \ref{orbit-theo}, since $P_{\ol{\alpha}}$ is the stabilizer subgroup in $H$ of the $H$-orbit corresponding to $\ol{\alpha}$.  Hence
 $T_p \mathcal{H}_{\ol{\alpha}}\stackrel{d\pi}{\simeq} T_{\pi(p)} M_{\ol{\alpha}} = T_{\pi(p)} M$, since $M_{\ol{\alpha}}$ is open.
This together with the
 $P$-equivariance of the curvature form of a Cartan connection, implies that the curvature form of the canonical Cartan connection $\omega^{nc}$ must vanish on $\mathcal{G}_{\vert M_0}$, which means it must vanish identically by continuity. Therefore, the conformal manifold $(M,[g])$ is locally conformally flat, and in particular its conformal holonomy $\mathrm{Hol}(M,[g])$ must be a discrete subgroup of $\SSO(3,2)$.

 \bbem
 This proof does not immediately generalize to $\SSO(p,q)$ with  $p+q=n>3$.  We have seen in Theorem \ref{orbit-theo}  that for $n>3$ the union $M_0$ of maximal orbits $M_{\ol{\alpha}}$ is still dense in $M$,  but  the orbits themselves, being of codimension $n-3$, are not open anymore. Hence, we obtain a foliation of $M_0$ into initial submanifolds in directions of which the curvature $\Omega^{nc}$ of the normal conformal Cartan connection $\w^{nc}$ on $\pi:\cal G|_{M_0}\to M_0$ vanishes. More precisely, we obtain that
 \[\Omega^{nc}_p(U,V)=0,\]
 for all $p\in \cal G$ with $\pi(p)\in M_0$ and $U,V\in T_p\cal G$ such that $d\pi_p(U)$ and $\d\pi_p(V)$ are tangent to an orbit $M_{\ol{\alpha}}$. This translates to the property that the Weyl tensor $W$ vanishes tangential to $M_{\ol{\alpha}}$, i.e.,
 \begin{equation}\label{weyl}W(U,V,X,Y)=0,
 \end{equation}
 for all $X,Y\in TM$ and $U,V\in TM_{\ol{\alpha}}$. Clearly, in general this does not force $W$ to vanish in all directions. However, for $n=4$ one can show that it does, and we sketch the argument here, saving the full details for future work in which we also plan to study the induced structures on the orbits of maximal dimension for arbitrary $n$ and the cases $H=\PSU(p,q)$ and $H=\PSp(p,q)$. For $n=4$ and $H=\SSO(3,1)\subset \SSO(6,3)$ or $\SSO(2,2)\subset \SSO(5,4)$, the maximal orbits are of codimension one. Given equation \eqref{weyl} and the first  Bianchi-identity for $W$,  for showing that $W\equiv 0$  it suffices to verify that $W(U,X,V,X)=0$ for $U,V\in TM_{\ol{\alpha}}$ and $X$ transversal to $M_{\ol{\alpha}}$. But  this follows from relation \eqref{weyl} and from $W$ having zero trace, {\em provided that the conformal metric   remains non-degenerate when restricted to  $TM_{\ol{\alpha}}$}. This is equivalent to the property that the  maximal $H$-orbits in the
M\"obius sphere are non-degenerate for the flat conformal metric, which is verified by a straightforward computation.
\ebem

\subsection{Proof of Theorem \ref{sutheo}} Let $(M,[g])$ be a conformal manifold of signature $(3,3)$, and follow the notational conventions of Section \ref{sec-background3} (so $G:=\SSO(4,4)$, $P \subset G$ is the parabolic subgroup stabilizing a null line in $\R^{4,4}$, etc.). If we denote by $H \subset G$ the image of the isotropy representation $\Ad_{\SSL_3\C}: \SSU(2,1) \rightarrow \SSO(4,4)$, then $H \simeq \PSU(2,1)$ and, by Proposition \ref{unitary prop}, we know that the open $H$-orbits are dense in $G/P$ (which is the double cover of $\SS^{3,3}$). Moreover, the stabilizer subgroup of an open $H$-orbit is given by $B \simeq \mathrm{U}(1) \times \mathrm{O}(1,1)$ and Proposition \ref{unitary prop} gives the explicit representation of $\b \subset \h$.

Now suppose $\mathrm{Hol}(M,[g]) \subseteq H \subset \SSO(4,4)$. Then applying Theorem \ref{cgh theorem} as in the proof of Theorem \ref{gensotheo} above, we see that there is a dense subset $M_0 \subset M$, consisting of the curved orbits corresponding to the open $H$-orbits in $G/P$,  and a canonical Cartan geometry $(\mathcal{H}_0 \rightarrow M_0,\eta_0)$ which reduces the canonical conformal Cartan geometry $(\mathcal{G} \rightarrow M,\omega^{nc})$ to type $(H,B)$, and, by Theorem \ref{Cartan invariant conformal holonomy}, induces a canonical metric $g_0 \in [g_{\vert M_0}]$ and a metric connection $\nabla^0$ with totally skew-symmetric, $\nabla^0$-parallel torsion $T^0$. Indeed, it is a straightforward matter to verify that $(\n,K_{\g} \vert_{\n})$ and $(\widehat{\n}, K_{\g} \vert_{\widehat{\n}})$ are homothetic as required by Theorem 7 (cf. also the calculations for $\sp(2,1)$ in Appendix B). Now we claim that $g_0$ has a nearly para-K\"ahler structure with canonical connection $\nabla^0$. Recall that an {\em almost para-K\"ahler structure} is an endomorphism field $J$ on a manifold $M_0$ with metric $g_0$ of neutral signature, such that $J$ squares to the identity, has two eigen distributions of the same rank,  and such that $J^*g_0=-g_0$. A {\em nearly para-K\"ahler structure} is an almost para-K\"ahler structure such that $(\nabla_XJ)(X) =0$ for all $X\in TM$ and $\nabla$ the Levi-Civita connection of $g_0$. A nearly para K\"ahler structure is {\em of constant type $\Lambda$} if
\begin{equation}\label{constanttype}g_0\left( (\nabla_XJ)Y, (\nabla_XJ)Y\right)=\Lambda\left(g_0(X,X)g_0(Y,Y)-(g_0(X,Y))^2+(g_0(JX,Y))^2\right),\end{equation}
for a constant $\Lambda$.

We begin by fixing the basis $\{u_1,u_2,u_3\}$ of $\C^{2,1}$ as in Proposition \ref{unitary prop}, in which the Hermitian form  is given as
\[ \la. , .\ra = \begin{pmatrix} 0& 0 & 1 \\
                       0 & 1 & 0 \\
                       1 & 0 & 0
                       \end{pmatrix}=:T.\]
We define the group 
$H:=\{A\in \SSL_3\C\mid \overline{A}^\top T A=T\}
$, which is conjugated to $\SSU(2,1)$, and let $\h$ be its Lie algebra. Then
$\sl_3\C=\h \+ \m$ with
\[\m:= \{X\in \sl_3\C\mid \overline{X}^\top T= T X\}.\]

                  Furthermore, we fix a matrix $S = \mathrm{diag}(\mu,\lambda,\ol{\mu}) \in \m$, with $\mu \in \C \backslash \R$ and $\lambda \in \R \backslash \{0\}$, which is null with respect to the Killing form $K=K_{\sl_3\C}$ of $\sl_3\C$ (which we scale to the trace form). Letting $\widehat{S} := S^*$  be the conjugate-transpose matrix, then $\widehat{S}$ is also $K$-null and we may re-scale if necessary to ensure that $K(S,\widehat{S})=1$.

Then the stabilizer corresponding to the (open) $H$-orbit of the null ray $\R_+ S$ in $G/P$ is given by
$$B = \left\{ b_{\varphi,r} := \begin{pmatrix} re^{\mathrm{i}\vf} & 0 & 0 \\
                       0 &  e^{-2\mathrm{i}\vf} & 0 \\
                       0 & 0 & r^{-1}e^{\mathrm{i}\vf}
                       \end{pmatrix} \mid \varphi \in \R, r \in \R \backslash \{0\} \right\},$$
with Lie algebra $\mathfrak{b}$ given as the diagonal matrices in $\su(2,1)$ which, in the above basis, has the form
                       $$\h = \left\{\begin{pmatrix} \beta - \mathrm{i}\alpha & y & i\delta \\
                       x & 2i\alpha & -\overline{y} \\
                       i\gamma & -\overline{x} & -\beta - \mathrm{i}\alpha
                       \end{pmatrix} \mid \alpha, \beta, \delta, \gamma \in \R, x,y \in \C \right\}.$$
By Proposition \ref{reductive stabilizer}, we have a ($K$-)naturally reductive decomposition $\h = \mathfrak{b} \dsum^\bot \mathfrak{n}$, with $\n = \b^{\perp} \subset \h$ given explicitly as
\[
 \mathfrak{n} = \left\{v(x,y,\gamma,\delta):=\begin{pmatrix} 0 & y & \mathrm{i}\delta\\
                       x & 0 & -\overline{y} \\
                       \mathrm{i}\gamma & -\overline{x} & 0
                       \end{pmatrix} \mid  x,y \in \C, \gamma, \delta \in \rr \right\}.
                       \]
The restriction of $K$ to $\n$ is given by the quadratic form \begin{equation}\label{null}
2(xy+ \overline{xy}-\gamma\delta),\end{equation} which has
signature $(3,3)$. In order to define the para-K\"ahler structure on $M_0$ we decompose $\mathfrak{n}$ further into
\be
\mathfrak{n}_+ & := &\left\{v_+(x,\delta):=\begin{pmatrix} 0 & 0& \mathrm{i}\delta \\
                       x & 0 & 0 \\
                       0 & -\ol{x} & 0
                       \end{pmatrix} \mid x \in \C, \delta\in \rr\right\},\\
\mathfrak{n}_- &:=&\left \{v_-(y,\gamma):=\begin{pmatrix} 0 & y & 0 \\
                       0 & 0 & -\ol{y} \\
                       \mathrm{i}\gamma & 0 & 0
                       \end{pmatrix} \mid y \in \C, \gamma\in \rr\right\}.
                       \ee
In particular, note that with this notation we can rewrite (\ref{null}) as
\begin{align}
K(v_-(x,\delta),v_+(y,\gamma)) = K(v_+(x,\gamma),v_-(y,\delta)) = xy + \overline{xy} - \gamma \delta. \label{killing form expression}
\end{align}

\blem\label{lemman}
                       $\n_\pm$ are
totally null with respect to the Killing form $K$ of $\sl_3(\C)$, invariant under $\Ad_H(B)$ and satisfy
\begin{equation}\label{lemmaneq}
\left[ \n_+,\n_-\right]\subset \b,\ \ \ [\n_\pm,\n_\pm]\subset \n_\mp.\end{equation}
\elem
\bprf
It follows from \eqref{null}
that $\n_\pm$ are totally null.
A straightforward check that
$\Ad(b_{\vf,r})$
sends $ v(x,y,\gamma,\delta)$ to $v(\frac{1}{r}\e^{-3\mathrm{i}\vf}x, r\e^{-3\mathrm{i}\vf}y,\frac{1}{r^2} \gamma,r^2\delta)$ implies the $\Ad_H(B)$-invariance.
The remaining properties are computed straightforwardly:
\be
\left[ v_+(x,\delta),v_-(y,\gamma)\right]&=&
\begin{pmatrix} -\gamma\delta-xy & 0 & 0 \\
                       0 & xy-\ol{xy} & 0 \\
                       0 & 0 & \gamma\delta+\ol{xy}
                       \end{pmatrix} \ \in \ \b
                       \\[2mm]
\left[ v_+(x,\gamma),v_+(y,\delta)\right]&=&
v_-(\mathrm{i}(\delta\ol{x}- \gamma\ol{y}), x\ol{y}-\ol{x}y) \in \n_-
                       \\[2mm]
\left[ v_-(x,\gamma),v_-(y,\delta)\right]&=&
v_+(\mathrm{i}(\delta\ol{x}- \gamma\ol{y}), x\ol{y}-\ol{x}y) \in \n_+.
                       \ee
\eprf
\blem The splitting $\n=\n_+\+\n_-$ defines an almost para-K\"ahler structure $J$ on $(M_0,g_0)$.
\elem
\bprf
Using the isomorphisms $\psi_u: T_xM_0 \rightarrow \h/\b \simeq \n$, cf. (\ref{Cartan frames}), we can define a splitting of $TM_0$ into two null distributions $T^\pm$ via $\psi_u(T_{x}^\pm)=\n_\pm$. The $\Ad_H(B)$-invariance then ensures that this is independent of the chosen $u \in \cal H_x$. Now setting $J|_{T^\pm}=\pm \mathrm{Id}_{T^\pm}$ defines an almost para-K\"ahler structure with respect to the metric $g_0$.
\eprf
Now we will show that this almost para-K\"ahler structure is in fact nearly para-K\"ahler. As noted above, $(M_0,g_0)$ has a natural metric connection $\nabla^0$ with totally skew-symmetric, $\nabla^0$-parallel torsion $T^0$. Moreover, the torsion $T^0$ is given by
\begin{equation}\label{torsion}
\psi_u(T^0_{x}(X,Y))= -\left[ \psi_u(X),\psi_u(Y) \right]_{\n}.
\end{equation}
Since $\n_\pm$ are $\Ad_H(B)$-invariant and thus invariant under the holonomy of $\nabla^0$, the almost para-K\"ahler structure $J$ is parallel with respect to $\nabla^0$, i.e., $\nabla^0$ is the canonical connection for $J$. Note that the almost para-K\"ahler structure $J$ can be viewed as induced by the stable three-form that is defined by the torsion of the connection $\nabla^0$ (for stable forms and para-K\"ahler structures see for example \cite{clss-09}).
Now we are ready to show
\blem
The almost para-K\"ahler structure $(M_0,g_0,J)$ is nearly para-K\"ahler of constant type $\Lambda= \einhalb$.
\elem
\bprf
Denote by $\nabla$ the Levi-Civita connection of $g_0$. As $\nabla^0$ is a metric connection, we have that
\[\nabla^0-\nabla=\einhalb T^0,\]
cf. e.g. \cite[Corollary 2.1]{agricola-srni}. On the other hand, since $\n_\pm$ are $\Ad_H(B)$-invariant, the almost para-complex structure $J$ is parallel with respect to $\nabla^0$,
\[\nabla^0J=0.\]
Hence we get
\begin{equation}\label{nablaj}
(\nabla_XJ)(Y)=\einhalb\left(J(T^0(X,Y))- T^0(X,J(Y)\right).\end{equation}
Hence, for $Y=X$, the definition of $J$ implies
\[
(\nabla_XJ)(X)=-\einhalb T^0(X,J(X))= T^0(X_+,X_-) ,\]
where $X=X_++X_-$ with $\psi_u(X_\pm)\in \n_\pm$. Formula \eqref{torsion} for the torsion and $[\n_+,\n_-]\subset \b$ imply that $T^0(X_+,X_-)=0$ and thus, that $J$ is a nearly para-K\"ahler structure.

In order to verify that this nearly para-K\"ahler structure is of constant type we have to show that
\eqref{constanttype}
holds with $\Lambda=\einhalb$. (Recall that $K$ is scaled to the trace form.)

We set $X=X_++X_-$ and $Y=Y_++Y_-$ and compute, using that $\n_\pm$ are totally null, that the right-hand-side in \eqref{constanttype} is equal to
\begin{equation}\label{rhs} 4\Lambda \left( g_0 (X_+,X_-)g_0 (Y_+,Y_-)- \left(g_0 (X_+,Y_-)\right)^2 + \left(g_0 (X_-,Y_+)\right)^2\right).\end{equation}
In order to compute the left-hand-side in \eqref{constanttype} we write $T^0(X,Y)=T_+^0(X,Y)+ T^0_-(X,Y)$.
From equation \eqref{torsion} and \eqref{lemmaneq} in  Lemma \ref{lemman} we get
\begin{equation}\label{torsion+-}
\psi\left( T^0_\pm (X,Y)\right)\ =\ -\left[ \psi (X_\mp), \psi (Y_\mp)\right].
\end{equation}
Based on
\eqref{nablaj}, this, together with \eqref{torsion}  implies
\be
\lefteqn{
\| (\nabla_XJ)(Y)\|^2}\\
&=& g_0\left( (\nabla_XJ)(Y), (\nabla_XJ)(Y)\right)\\
&=&
\frac{1}{4}\left( \| J(T^0(X,Y)))\|^2 -2 g_0( J(T^0(X,Y)), T^0(X,JY ) + \| T^0(X,JY) \|^2 \right)
\\
&=&
-\einhalb\left(
g_0(T^0_+(X,Y),T^0_-(X,Y) - g_0(T^0_+(X,JY),T^0_-(X,JY) \right)
\\
&&
- \einhalb\left( g_0(T^0_+(X,Y),T^0_-(X,JY) - g_0(T^0_+(X,JY),T^0_-(X,Y) \right)
\\
&=&
-2K \left(\left[ \psi(X_-),\psi(Y_-)\right], \left[\psi( X_+),\psi(Y_+)\right]\right).
\ee
By using \eqref{killing form expression} and \eqref{lemmaneq} we can compare this to \eqref{rhs} and get that $\Lambda=\einhalb$. Alternatively, note that by Theorem \ref{IZtheo} cited below, any six-dimensional nearly para-K\"ahler manifold is automatically of constant type $\Lambda$ for some $\Lambda \in \R$, and so it suffices to verify the constant $\Lambda$ in \eqref{constanttype} by computing both sides of the equation for simple choices of $X$ and $Y$ for which \eqref{rhs} is non-zero, e.g. for $\psi(X) = v_+(1,0) + v_-(1,0)$ and $\psi(Y) = v_+(0,1) + v_-(0,1)$.
\eprf
The proof of Theorem \ref{sutheo} now follows from
\btheo[Ivanov \& Zamkovoy \cite{ivanov-zamkovoy05}] \label{IZtheo}
A six-dimensional nearly para-K\"ahler manifold is of constant type $\Lambda$ and Einstein with Einstein constant $5 \Lambda$.
\etheo
This implies that on the open and dense submanifold $M_0$ we have found an Einstein metric $g_0$ in the conformal class $[g]$, due to Theorem \ref{Cartan invariant conformal holonomy}, with positive Einstein constant $\frac{5}{2}$, which gives a parallel section of the tractor bundle over $M_0$ and  forces the normal conformal holonomy of $[g]$ over $M_0$ to be contained in the stabilizer of a time-like vector in $H \simeq \PSU(2,1) \subset \SSO(4,4)$.

\subsection{Proof of Theorem \ref{sptheo}} Let $(M^{5,7},[g])$ be a conformal manifold of signature $(5,7)$, and denote now by $H \subset \SSO(6,8)$ the image of the isotropy representation $\Ad_{\SSL_3\Q}: \SSp(2,1) \rightarrow \SSO(6,8)$. Then $H \simeq \PSp(2,1)=\PSp(2,1)/\{\pm\1\}$, and if $(M,[g])$ satisfies $\mathrm{Hol}(M,[g]) \subseteq H$, we have by the same arguments as in the proofs of Theorems \ref{gensotheo} and \ref{sutheo} a dense subset $M_0 \subset M$ and a canonical Cartan geometry $(\mathcal{H}_0 \rightarrow M_0,\eta_0)$ of type $(H,B)$. Due to Proposition \ref{symplectic prop}, $B $ is doubly covered by $\SSp(1) \times \SSL_2\C$. This is a reduction of the canonical (normal) conformal Cartan geometry $(\mathcal{G} \rightarrow M,\omega^{nc})$ of $(M,[g])$, and Theorem~\ref{Cartan invariant conformal holonomy} applies to show that $(\mathcal{H}_0,\eta_0)$ induces a canonical metric $g_0 \in [g\vert_{M_0}]$ and a metric connection $\nabla^0$ with skew-symmetric, parallel torsion $T^0$. That the canonical metric $g_0$ is in the conformal class of $g$ follows from the last part of Theorem~\ref{Cartan invariant conformal holonomy} and part~(a) of Lemma \ref{tech lemma} in Appendix~\ref{techlemmaappendix}.

We claim that the Ricci tensor of $g_0$ is given by the Ricci tensor of the naturally reductive pseudo-Riemannian space $\PSp(2,1)/B$. In particular, since $\PSp(2,1)/B$ is Einstein (cf. Appendix A), Theorem \ref{sptheo} follows from

\bs \label{Ricci prop} Let $K$ denote the naturally reductive metric for the decomposition $\sp(2,1) = \b \dsum \n$ given by Proposition \ref{reductive stabilizer}. The canonical metric $g_0$ as above has Ricci tensor $\mathrm{Ric}^{g_0}$ which is related to the Ricci tensor $\mathrm{Ric}^{K}$ of the naturally reductive homogeneous pseudo-Riemannian space $\PSp(2,1)/B$ by:
\begin{align}
\mathrm{Ric}^{g_0}(X,Y) = \mathrm{Ric}^{K}({\psi}_u(X),{\psi}_u(Y)), \label{Ricci formula}
\end{align}
for any $x \in M_0$, $X, Y \in T_xM$ and $u \in \mathcal{H}_x$, where $\psi_u$ is the map as given in (\ref {Cartan frames}) by the Cartan geometry $(\mathcal{H}_0,\eta_0)$.
\es

\begin{proof}  By part (b) of Lemma \ref{tech lemma} we have a $K$-orthonormal basis $\{e_1,\ldots,e_{12}\}$ of $\n$, and bases $\{E_1,\ldots,E_{12}\}$ of $\p_- \subset \so(6,8)$ and $\{E^1,\ldots,E^{12}\}$ of $\p_+ \subset \so(6,8)$, which are dual with respect to the Killing form $K_{\so(6,8)}$ of $\so(6,8)$, i.e. $K_{\so}(E_i,E^j) = \delta_{ij}$, and are related to the basis of $\n \subset \g$ by: $$E_i = ce_i + A_i \text{ and } E^i = \varepsilon_ice_i + A^i,$$ for $0 \neq c \in \R$, $A_i \in \p, A^i \in \widehat{\p}$, and $\varepsilon_i = K(e_i,e_i) = \pm 1$ (for notational simplicity, we identify $e_i$ and other elements of $\sp(2,1)$ with their image $\rho(e_i) \in \so(6,8)$, as this seems unlikely to cause confusion in the present context). Now, since the curvature form of the Cartan connection $\eta_0$ is just given by the restriction of the curvature form of $\omega^{nc}$ to the sub-bundle $\mathcal{H}_0 \subset \mathcal{G}$ over $M_0$, we can also identify their curvature functions $$\kappa = \kappa^{\eta_0} = \kappa^{nc} \vert_{\mathcal{H}_0}: \mathcal{H}_0 \rightarrow \Lambda^2\n^* \otimes \h \subset \Lambda^2(\so(6,8)/\p)^* \otimes \so(6,8),$$ (see Section \ref{sec-background}). Recall that the canonical conformal Cartan connection $\omega^{nc}$ satisfies the normality condition (\ref{normality}). We will now use the bases $\{e_i\}$, $\{E_i\}$ and $\{E^i\}$ to translate the normality of $\omega^{nc}$ into a geometric condition for the reductive Cartan connection $\eta_0$. For the following calculation, note in particular that $\kappa(u) \in \Lambda^2\n^* \otimes \b$ since the torsion of $\eta_0$ (and of $\omega^{nc}$) vanishes, and that $\kappa(u)(Y,.) = 0$ for all $Y \in \p$. Then we get, for any $X \in \n$:
\begin{align*}
0 = (\partial^* \circ \kappa)(u)(X) &:= \sum_{i=1}^{12}[\kappa(u)(X,E_i),E^i] = \sum_{i=1}^{12}[\kappa(u)(X,ce_i+A_i),\varepsilon_ice_i+A^i] \\
&= c^2\sum_{i=1}^{12}\varepsilon_i[\kappa(u)(X,e_i),e_i] + c\sum_{i=1}^{12}[\kappa(u)(X,e_i),A^i].
\end{align*}
But now we note, since $\kappa(u)(X,e_i) \in \b$ for all $i=1,\ldots,12$, that all terms of the first sum in the final line must lie in $\n$ (since $[\b,\n] \subset \n$ by reductivity).  Now the results in Appendix \ref{techlemmaappendix}
show on one hand that $[\b,\widehat{\p}] \subset \widehat{\p}$ as a result of $\b \subset \p_0$, which implies that 
   all terms of the second sum must lie in $\widehat{\p}$, and on the other that  
   $\n \cap \widehat{\p} = \{0\}$, as a result of $\b = \sp(2,1) \cap \widehat{\p}$. Hence, 
    each sum must vanish separately.

Now, the vanishing of the sum
\begin{align}
\sum_{i=1}^{12}\varepsilon_i[\kappa(u)(X,e_i),e_i], \label{reduced sum}
\end{align}
for $X \in \n$ and $e_1,\ldots,e_{12}$ an orthonormal basis of $\n$ has a very natural geometric meaning for the reductive Cartan geometry $(\mathcal{H}_0,\eta_0)$, namely we claim it means that the Ricci tensor of the covariant derivative $\nabla^0$ which $\eta_0$ induces on $M_0$ is equal to the Ricci tensor of the natural covariant derivative on $\PSp(2,1)/B$. To see this, note that the curvature tensor $R^0$ of $\nabla^0$ satisfies the formula:
$$\psi_u \circ R^{0}(X,Y) \circ \psi_u^{-1} = \Omega^{\eta_0}(\widehat{X},\widehat{Y}) - [\psi_u(X),\psi_u(Y)]_{\b}$$ for any $X, Y \in T_xM_0$, $u \in \mathcal{H}_x$ and $\widehat{X}, \widehat{Y} \in T_u \mathcal{H}_0$ projecting to $X, Y$. Note that the second term on the right-hand-side equals $R^{\n}(\psi_u(X),\psi_u(Y))$, where $R^{\n}$ is the curvature of the natural covariant derivative $\nabla^{\n}$ of $\PSp(2,1)/B$ which has torsion $T^{\n}({\psi}_u(X),{\psi}_u(Y)) = -[{\psi}_u(X),{\psi}_u(Y)]_{\n}$. Hence it follows, using the relation between $\kappa$ and $\Omega^{\eta_0}$, that the vanishing of the sum (\ref{reduced sum}) implies that the Ricci tensor $\mathrm{Ric}^{0}$ of $\nabla^0$ satisfies
\begin{align}
\mathrm{Ric}^{0}(X,Y) = \mathrm{Ric}^{\n}(\psi_u(X),\psi_u(Y)), \label{first ricci formula}
\end{align}
where $\mathrm{Ric}^{\n}$ is the Ricci tensor of $\nabla^{\n}$. Now the result relating the Ricci tensors
$\mathrm{Ric}^{g_0}$ and $\mathrm{Ric}^{K}$ of the Levi-Civita connections of $g_0$ and $K$ follows, since the torsions $T^{0}$ and $T^{\n}$ of $\nabla^0$ and $\nabla^{\n}$, respectively, are both totally skew-symmetric and parallel. Hence, from the formula in \cite[Proposition 3.1]{ivanov-papa01} that relates $\mathrm{Ric}^{g_0}$ and $\mathrm{Ric}^0$ (see also \cite[Theorem A.1]{agricola-srni}), and using that $\nabla^0T^0=0$ implies that the $\nabla^0$-divergence of $T^0$ vanishes, we get:
$$\mathrm{Ric}^{g_0}(X,Y) = \mathrm{Ric}^{0}(X,Y) - \frac{1}{4}\sum_{i=1}^{12}\varepsilon_i g_0(T^{0}(u_i,X),T^{0}(u_i,Y)),$$ for any $g_0$-orthonormal basis $\{u_i\}$ of $T_pM_0$. An analogous formula relates $\mathrm{Ric}^{K}$ and $\mathrm{Ric}^{\n}$. Now the claimed identity (\ref{Ricci formula}) follows from the relations $$\psi_u(T^{\nabla}(X,Y)) = T^{\n}({\psi}_u(X),{\psi}_u(Y)) \,\,\, \mathrm{and} \,\,\, g_0(X,Y) = K({\psi}_u(X),{\psi}_u(Y)),$$ completing the proof. \end{proof}

\section{Ambient extension and a general non-existence result}
\label{fg}

In this section, we prove Theorem \ref{ambient holonomy extension}. The proof relies on a result concerning the so-called \emph{ambient extension} of parallel tractors for a conformal manifold $(M,[g])$ to parallel tensors of the Fefferman-Graham ambient space $(\widetilde{M},\widetilde{g})$ of $(M,[g])$, so we briefly review the necessary facts.

\subsection{The Fefferman-Graham ambient space for odd-dimensional real-analytic conformal manifolds}

The following summary is based on  \cite{fefferman/graham85} and \cite{fefferman-graham07}, cf. also Section \ref{sec-background} of \cite{graham-willse11}: For a conformal pseudo-Riemannian manifold $(M,[g])$ of signature $(p,q)$, we have the principal $\R_+$-bundle $\pi: \mathcal{Q} \rightarrow M$ defined by $\mathcal{Q} = \{ (p,g_p) : p \in M, g \in [g] \}$, where $\pi$ is the canonical projection and the right $\R_+$-action is given by dilation $\delta_s(p,g_p) = (p,s^2 g_p)$. The tautological tensor $\mathbf{g}$ is a degenerate symmetric bilinear form on $\mathcal{Q}$ defined by $\mathbf{g}_{(g_p,p)}(U,V) := g_p(\pi_*(U),\pi_*(V))$. Extend the $\R_+$-action on $\mathcal{Q}$ to $\mathcal{Q} \times \R$ and let $\widetilde{M} \subset \mathcal{Q} \times \R$ be an $\R_+$-invariant open subset containing the inclusion $\iota(\mathcal{Q}) = \mathcal{Q} \times \{ 0 \}$. Then a \emph{pre-ambient metric} for $(M,[g])$ is given by some smooth pseudo-Riemannian metric $\widetilde{g}$ on $\widetilde{M}$ of signature $(p+1,q+1)$, which satisfies (i) $\delta_s^*\widetilde{g} = s^2 \widetilde{g}$ for $s \in \R_+$; and (ii) $\iota^*\widetilde{g} = \mathbf{g}$. A pre-ambient metric is called \emph{straight} if the flow by dilation, $s \mapsto \delta_s(p)$, is a geodesic with respect to $\widetilde{g}$ for all $p \in \widetilde{M}$ (equivalently, if the fundamental vector field of the dilation action, $T = \frac{d}{ds}\delta_s \vert_{s=1}$, satisfies $\widetilde{\nabla} T = \mathrm{Id}$ for the Levi-Civita connection of $\widetilde{g}$, cf.  \cite[Propositions 2.4 and 3.4]{fefferman-graham07}). An \emph{ambient metric} for $(M,[g])$ is then defined to be a pre-ambient metric with Ricci tensor vanishing to certain orders (with respect to the $\R$-component of $\widetilde{M}$) depending on whether the dimension $n = p+q$ is even or odd, and the pair $(\widetilde{M},\widetilde{g})$ is called an \emph{ambient space} for $(M,[g])$. Here we do not re-state the conditions in full generality, but will only consider the odd-dimensional case where $(M,[g])$ is real-analytic (i.e., some $g \in [g]$ is real-analytic), where the questions of existence and uniqueness are simplified and the following fundamental result holds (cf. \cite{fefferman/graham85,kichenassamy04,fefferman-graham07}):

\btheo[Fefferman \& Graham \cite{fefferman/graham85,fefferman-graham07}]  \label{Fefferman Graham thm}
Let $(M,[g])$ be a real-analytic conformal manifold of odd dimension $n=p+q > 1$. Then there exists an ambient space $(\widetilde{M},\widetilde{g})$ for $(M,[g])$ with real-analytic Ricci-flat metric $\widetilde{g}$. The ambient space is unique modulo diffeomorphisms that restrict to the identity along $\iota(\mathcal{Q}) \subset \widetilde{M}$ and commute with the $\R_+$-action.
\etheo

Starting from the Fefferman-Graham ambient space $(\widetilde{M},\widetilde{g})$, it was shown by \v{C}ap and Gover in \cite{cap-gover03} that the tangent bundle $T\widetilde{M}$ and Levi-Civita connection $\widetilde{\nabla}$ of $(\widetilde{M},\widetilde{g})$ induce the standard conformal tractor bundle $\mathcal{T}$ and normal tractor connection $\nabla^{\mathcal{T}}$ of $(M,[g])$ by identifying $$\mathcal{T}_x \cong \{ U \in \Gamma(T\widetilde{M}_{\vert \mathcal{Q}_x}) : [T,U] = -U \}$$ and showing that the properties of $\widetilde{g}$ imply that $\widetilde{\nabla}$ descends to a well-defined linear connection $\Gamma(\mathcal{T}) \rightarrow \Gamma(T^*M \otimes \mathcal{T})$ satisfying the normalization condition which uniquely determines $\nabla^{\mathcal{T}}$ (actually, in \cite{cap-gover03} it is shown that this still holds under a weakening of the conditions on the ambient space $(\widetilde{M},\widetilde{g})$).

Note that from the straight-ness of the ambient metric $\widetilde{g}$, i.e. the property $\widetilde{\nabla} T = \mathrm{Id}$, it is an easy consequence to see that a $\widetilde{\nabla}$-parallel vector field on $\widetilde{M}$ automatically determines a section of the tractor bundle $\mathcal{T}$, which by the above is parallel with respect to the normal conformal tractor connection. Similarly, identifying tensor powers of the standard conformal tractor bundle as $$\bigotimes^k \mathcal{T}_x^* \cong \{ \Upsilon \in \Gamma(\bigotimes^k T^*\widetilde{M}_{\vert \mathcal{Q}_x}) : \mathcal{L}_T \Upsilon = k\Upsilon \},$$ one sees that $\widetilde{\nabla}$-parallel tensors on $\widetilde{M}$ restrict to parallel tractors. Moreover, the conformal holonomy group of $(M,[g])$ is contained in the (pseudo-Riemannian) holonomy group of $(\widetilde{M},\widetilde{g})$,
\begin{align}
\mathrm{Hol}(M,[g]) \subseteq \mathrm{Hol}(\widetilde{M},\widetilde{g}), \label{tractor vs ambient holonomy}
\end{align}
(cf. \cite[Proposition 6.2]{armstrong-leistner07}  where a more general result is shown, which recovers this inclusion by noting that the Levi-Civita connection $\widetilde{\nabla}$ of the ambient space automatically satisfies the assumptions).

On the other hand, given a tractor $\Upsilon \in \Gamma(\bigotimes^k \mathcal{T}^*)$, an \emph{ambient extension} is defined to be a tensor $\widetilde{\Upsilon} \in \Gamma(\bigotimes^k T^*\widetilde{M})$ which satisfies $\delta_s^*\widetilde{\Upsilon} = s^k\widetilde{\Upsilon}$ and $\widetilde{\Upsilon}\vert_{\mathcal{Q}} = \Upsilon$ (cf. \cite[Section \ref{sec-orbit}]{graham-willse11}), and one can ask whether a $\nabla^{\bigotimes^k\mathcal{T}^*}$-parallel tractor $\Upsilon$ has an ambient extension which is $\widetilde{\nabla}$-parallel to some order. This problem was studied in \cite{graham-willse11}, where it was proved (again, we cite only the result for $n$ odd and $(M,[g])$ real-analytic, noting that results were also obtained under weaker assumptions):

\btheo[Graham \& Willse \cite{graham-willse11}] \label{Graham Willse thm}
Let $(M,[g])$ be a real-analytic conformal manifold of odd dimension $n > 1$, and let $\widetilde{g}$ be a real-analytic Ricci-flat ambient metric for $(M,[g])$. If $\Upsilon \in \Gamma(\bigotimes^k \mathcal{T}^*)$ is parallel with respect to the normal conformal tractor connection, then $\Upsilon$ has a real-analytic ambient extension $\widetilde{\Upsilon}$ satisfying $\widetilde{\nabla} \widetilde{\Upsilon} = 0$ in a neighborhood of $\mathcal{Q} \times \{ 0 \}$ in $\widetilde{M}$.
\etheo

\subsection{Proof of Theorem \ref{ambient holonomy extension}} Let $(M,[g])$ be an odd-dimensional, real-analytic conformal manifold, that is the underlying manifold $M$ is real-analytic and there is a metric $g \in [g]$ whose coefficients with respect to any real-analytic local chart of $M$ are real-analytic. By Theorem \ref{Fefferman Graham thm}, a Ricci-flat, real-analytic ambient metric $\tilde{g}$ exists on some ambient space $\widetilde{M} \approx \R_+ \times M \times \R$. Let $\Upsilon \in \Gamma(\bigotimes^k \mathcal{T}^*)$ be a parallel tractor determining the conformal holonomy reduction $\mathrm{Hol}(M,[g]) \subseteq H$, where $H$ is the identity component of the stabilizer in $\O(p+1,q+1)$ of some vector in $\bigotimes^k (\R^{p+1,q+1})^*$. Then by Theorem \ref{Graham Willse thm}, $\Upsilon$ has a parallel ambient extension to $\widetilde{M}$, and therefore the pseudo-Riemannian holonomy of the ambient space is reduced to the isotropy subgroup $H$, $\mathrm{Hol}(\widetilde{M},\widetilde{g}) \subseteq H$.

Now suppose, in contradiction to the statement of Theorem \ref{ambient holonomy extension}, that $\mathrm{Hol}(M,[g]) = H$. Then by the inclusion $\mathrm{Hol}(M,[g]) \subseteq \mathrm{Hol}(\widetilde{M},\widetilde{g})$, cf. (\ref{tractor vs ambient holonomy}), we must have $\mathrm{Hol}(\widetilde{M},\widetilde{g}) = H$. Then, when consulting Berger's list of irreducible, non-symmetric pseudo-Riemannian holonomy groups, we find that either $H=\SSO^0(p+1,q+1)$, which is the stabilizer of the curvature tensor of a space of constant curvature, or $H=\G_{2(2)}$, which is the stabilizer of a stable $3$-form on $\R^{3,4}$, or that $(\tem, \tg)$ is a locally symmetric space and thus locally isometric to an irreducible, non-flat pseudo-Riemannian symmetric space. Hence, by Proposition \ref{Ricci nonflat}, the Ricci tensor $\mathrm{Ric}_{\widetilde{g}}$ is non-zero, a contradiction to the defining properties of $\widetilde{g}$.\hfill $\Box$

\begin{appendix}

\section{The naturally reductive space $\PSp(2,1)/B$}
\label{appA}
In this and the following appendix we will work in a basis that was obtained in the proof of Proposition~\ref{symplectic prop}, and for which the Hermitian form on $\Q^3$ is of the form 
\[T:=\begin{pmatrix} 0& 0 & 1 \\
                       0 & 1 & 0 \\
                       1 & 0 & 0
                       \end{pmatrix}.\]
Hence, our conventions here are slightly different from Section \ref{secgroupdefs}: We define
\[
\wt H:=
\left\{
\begin{pmatrix} A&-B\\\ol{B}&\ol{A}\end{pmatrix}
\mid \ol{A}^\top T A+B^\top T\ol{B}=T,\ B^\top T  \ol{A}- \ol{A}^\top T B=0\right\}\subset \SSU^*(6)\simeq\SSL_3\Q,
\]
 which is conjugated in $ \SSU^*(6)$ to $\SSp(2,1)$ as defined in Section \ref{sympsec}. Its Lie algebra is given as
 \[
 \h
 =
 \left\{\begin{pmatrix} X & -Y \\ \ol Y&\ol X\end{pmatrix}\in \sl_6\C \mid  X,Y\in \gl_n\C, \  \ol{X}^\top T+ T X=0, \ Y^\top T-T Y=0\right\}\simeq \sp(2,1),
\]
and we have $\su^*(6)=\h\+\m$ with $\Ad(\wt H)$-invariant
\[
 \m
 =
 \left\{\begin{pmatrix} X & -Y \\ \ol Y&\ol X\end{pmatrix}\in \sl_6\C \mid  X,Y\in \gl_n\C,  \  \ol{X}^\top T- T X=0,  \ Y^\top T+T Y=0\right\}.
\]
Now let $H:=\Ad (\wt H)\simeq \wt H/\{\pm \1\}\simeq \PSp(2,1)$  be the image of $\wt H$ of the adjoint action on $\m$, and $[S] \in \mathbb{S}(\mathfrak{m})_0$ a null line as in Proposition \ref{symplectic prop} spanned by
\begin{equation}\label{Seqn}S = \begin{pmatrix} S_0 & 0 \\ 0 & \overline{S_0} \end{pmatrix} \in \m,\end{equation}
 where $S_0 = \mathrm{diag}(\mu,-2\mathrm{Re}(\mu),\overline{\mu}) \in \gl_3\C$ for a suitable choice of complex number $\mu$ such that $S^2$ has no trace.
Consider the homogeneous space $H/B$ for $B \simeq (\SSp(1) \times \SSL_2\C)/\{\pm \1\}$ the stabilizer in $H$ of $\R S$. Let $\h = \b \dsum \n$ be the reductive decomposition given by applying Proposition \ref{reductive stabilizer}, and $K = K_{\su^*(6)} \vert_{\h}$ the corresponding naturally reductive metric. Then $K$ and $K_{\h}$ are both given, up to a multiple, by the trace form over $\C^6$. Explicitly, we have
\begin{equation}
\b = \left\{ B(z,ix,y_1,y_2,y_3) := \begin{pmatrix}X&-Y\\\ol{Y}&\ol{X}\end{pmatrix} \left| 
\begin{array}{l}
X = \mathrm{diag}(z,ix,-\ol{z}) , \ Y = \mathrm{adiag}(y_1,y_2,y_3),
\\
z,y_1,y_2,y_3 \in \C, x \in \R
\end{array} \right. \right\}
\label{bdef}
\end{equation}
and its complement is given as 
\begin{equation}
\label{ndef}
\n = \left\{ \begin{pmatrix}X&-Y\\\ol{Y}&\ol{X}\end{pmatrix} \left| X=
\begin{pmatrix}
0 & z_1 & ix_1\\
z_2 & 0 & -\ol{z_1} \\
ix_2 & -\ol{z_2} & 0
\end{pmatrix},
Y=
\begin{pmatrix}
y_1 & y_2 & 0\\
y_3 & 0 & y_2\\
0 & y_3 & y_1
\end{pmatrix},\ z_i,y_i \in \C, x_1,x_2 \in \R
\right.
\right\}.
\end{equation}
We further decompose $\b = \b_1 \dsum \b_2$ by letting
\begin{align}
\b_1 &= \{ B_1(ix,y) := B(0,ix,0,y,0) \} \simeq \sp(1); \label{b1 def}\\
\b_2 &= \{ B_2(z,y,w) := B(z,0,y,0,w) \} \simeq \sl_2\C. \label{b2 def}
\end{align}
And we split $\n = \n_1 \dsum \n_2$ as
\begin{align}
\n_1 &= \{ x_1 = x_2 = y_1 = 0 \} \simeq \{ (z_1,z_2,y_2,y_3)^{\top} \in \C^4 \}; \label{n1 def}\\
\n_2 &= \{ z_1=z_2=y_2=y_3 = 0 \} \simeq \{ (ix_1,ix_2,y_1)^{\top} \in \mathrm{Im}(\C)^2 \dsum \C \}. \label{n2 def}
\end{align}

Then a straightforward computation shows that $\b_1$ acts trivially on $\n_2$ via the adjoint action, while the action on $\n_1$ is given by
\begin{align}
\ad(B_1(ix,y)): \begin{pmatrix} z_1 \\ z_2 \\ y_2 \\ y_3 \end{pmatrix} \mapsto \begin{pmatrix} -ixz_1 + \ol{y}y_2 \\ ixz_2 - y\ol{y_3} \\ ixy_2 - yz_1 \\ ixy_3 + y\ol{z_2} \end{pmatrix}. \label{b1 action on n1}
\end{align}
Similarly, we calculate that $\b_2$ preserves the decomposition $\n = \n_1 \dsum \n_2$ under the adjoint action, and the action on the $\n_1$ and $\n_2$ summands is given, respectively, by
\begin{align}
\ad(B_2(z,y,w)): \begin{pmatrix} z_1 \\ z_2 \\ y_2 \\ y_3 \end{pmatrix} &\mapsto \begin{pmatrix} zz_1 - y\ol{y_3} \\ -zz_2 + \ol{w}y_2 \\ zy_2 - yz_2 \\ -\ol{z}y_3 + w\ol{z_1} \end{pmatrix}; \label{b2 action on n1} \\
\ad(B_2(z,y,w)): \begin{pmatrix} ix_1 \\ ix_2 \\ y_1 \end{pmatrix} &\mapsto \begin{pmatrix} i2(\mathrm{Re}(z)x_1 + \mathrm{Im}(\ol{y}y_1)) \\ -i2(\mathrm{Re}(z)x_2 - \mathrm{Im}(\ol{w}y_1)) \\ i2\mathrm{Im}(z)y_1 - wix_1 - yix_2 \end{pmatrix}. \label{b2 action on n2}
\end{align}

Fixing $K$ to be one-half the trace form over $\C^6$, one verifies directly that a $K$-orthonormal basis of $\b_1$ is given by:
\begin{align}
A_1 := B_1(i,0), \,\, A_2 := B_1(0,1), \,\, A_3 := B_1(0,i); \label{b1 basis}
\end{align}
and these satisfy $K(A_i,A_j) = \varepsilon_i\delta_{ij}$ for $\varepsilon_i=-1$, $i,j=1,2,3$. Also, a $K$-orthonormal basis of $\b_2$ is given by
\begin{align}
A_4 := \frac{1}{\sqrt{2}}B_2(i,0,0), \,\, A_5 &:= \frac{1}{\sqrt{2}}B_2(0,1,1), \,\, A_6 := \frac{1}{\sqrt{2}}B_2(0,i,i), \\
A_7 := \frac{1}{\sqrt{2}}B_2(1,0,0), \,\, A_8 &:= \frac{1}{\sqrt{2}}B_2(0,1,-1), \,\, A_9 := \frac{1}{\sqrt{2}}B_2(0,i,-i); \label{b2 basis}
\end{align}
and these satisfy $K(A_i,A_j) = \varepsilon_i\delta_{ij}$, where $\varepsilon_i = -1$ for $i=4,5,6$ and $\varepsilon_i = 1$ for $i=7,8,9$.

Hence the Casimir operator of the representation $\rho = \ad_{\sp(2,1)}: \b \rightarrow \gl(\n)$ with respect to $K$ is given, up to sign, using the above basis of $\b$, as:
\begin{align*}
\chi_{\rho,K} = \sum_{i=1}^9 \varepsilon_i \rho(A_i) \circ \rho(A_i).
\end{align*}
Now it is only mildly tedious, and perhaps even enjoyable, to calculate, using the definitions (\ref{b1 basis})-(\ref{b2 basis}) and the formulae (\ref{b1 action on n1})-(\ref{b2 action on n2}), the identity
\begin{align*}
\chi_{\rho,K} = 6\mathrm{Id}_{\n}. 
\end{align*}
No we apply the result by Wang and Ziller  \cite[page 569, (1.7) Corollary]{wang-ziller85},
that a naturally homogeneous metric is Einstein if and only if $\chi_{\rho,K}$ is a multiple of the identity. Hence, we obtain
 that the naturally reductive homogeneous metric, induced on $\PSp(2,1)/B$ by $K$, is indeed Einstein.

\section{Proof of technical Lemma \ref{tech lemma}}
\label{techlemmaappendix}
With  $\g:=\su^*(6) = \halg \dsum \malg \subset \sl_6\C$ the symmetric decomposition wit $\h$ and $\m$ defined in the previous appendix, let $\som$ denote the special orthogonal algebra of the restriction of the Killing form $K_{\galg}$ of $\galg$ to $\mathfrak{m}$, and denote by $\rho: \halg \rightarrow \som$ the isotropy representation, which is faithful. In Section 3, we identified a null vector $S \in \malg$ explicitly given  in \eqref{Seqn},
 and a Cartan involution $\theta$ of $\halg$ defined by $\theta: X \mapsto -\overline{X}^{\top}$, such that the stabilizer subalgebra $\b = \mf{stab}_{\h}(\R S)$ as in \eqref{bdef}  is $\theta$-invariant. This implies that $\b$ stabilises another  line spanned by  $\widehat{S} = \theta(S)$. We have that $\widehat{S} \in \malg$ is also a null vector and $K_{\g}(S,\widehat{S}) \neq 0$ (since $K_{\g}(S,\widehat{S}) = K_{\g}(S,\widehat{\theta}(S)) < 0$, we will also assume after rescaling if necessary, that $K_{\g}(S,\widehat{S}) = -1$). Therefore, $S$ and $\widehat{S}$ define a \emph{$\vert 1 \vert$-grading} of $\som$: Let $\p := \mf{stab}_{\so(\m)}(\R S)$, $\widehat{\p} := \mf{stab}_{\so(\m)}(\R \widehat{S})$ and define $\p_+ := \mathrm{Ker}(\ol{\ad}: \p \rightarrow \gl(\so(\m)/\p))$, $\p_- := \mathrm{Ker}(\ol{\ad}: \widehat{\p} \rightarrow \gl(\so(\m)/\widehat{\p}))$, and $\p_0 := \p \cap \widehat{\p}$. This determines a vector space decomposition
$$\so(\m) = \p_- \dsum \p_0 \dsum \p_+$$
which is a $\vert 1 \vert$-grading.
Then the $\theta$-invariance of $\b=\h\cap \p$ implies 
$\b\subset \p_0$ 
and $\b=\h\cap \widehat\p$.

We also have $\widehat{\n} := \mf{s}^{\perp} \subset \m$, for $\mf{s} = \mathrm{span}(S,\widehat{S})$, and we have a reductive decomposition $\h = \b \dsum \n$ where $\b = \mf{stab}_{\h}(\R S)$. We will prove:

\blem \label{tech lemma} \noindent (a) There is a linear isomorphism $\n \simeq \widehat{\n}$ which pulls back $K_{\g}\vert_{\widehat{\n}}$ to $K_{\g}\vert_{\n}$.

\noindent (b) We can find $K_{\g}$-orthonormal bases $\{e_1,\ldots,e_{12}\}$ of $\n$, and $K_{\so(\m)}$-dual bases $\{E_1,\ldots,E_{12}\}$ of $\p_-$ and $\{E^1,\ldots,E^{12}\}$ 
of $\p_+$, which are related, for some constant $c$, by
\begin{align}
\frac{1}{c}\rho(e_i) = E_i + A^i +\varepsilon_iE^i, \label{relate ei to Ei}
\end{align}
for $A^i \in \p_0$, where $\varepsilon_i := K_{\g}(e_i,e_i) = \pm 1$.
\elem

\begin{proof} For part (a),  using the form of $S \in \m$ in \eqref{Seqn} and  $\widehat{S} = -\overline{S}^{\top} = -\ol{S}$, one checks  that
\[
\widehat{\n} = \left\{ \begin{pmatrix}X&-Y\\\ol{Y}&\ol{X}\end{pmatrix} \left| X=
\begin{pmatrix}
0 & x_1 & \alpha_1\\
x_2 & 0 & \ol{x_1} \\
\alpha_2 & \ol{x_2} & 0
\end{pmatrix},
Y=
\begin{pmatrix}
y_1 & y_2 & 0\\
y_3 & 0 & -y_2\\
0 & -y_3 & -y_1
\end{pmatrix}
\right.
\right\},
\]
for $x_1,x_2,y_1,y_2,y_3 \in \C, \alpha_1,\alpha_2 \in \R$. We fix $K_{\g}$ to be one-half the trace form (over $\C^6$), and note
\begin{align}
K_{\g} \left( \begin{pmatrix}X&-Y\\\ol{Y}&\ol{X}\end{pmatrix}, \begin{pmatrix}V&-W\\\ol{W}&\ol{V}\end{pmatrix} \right) = \mathrm{Re}(\mathrm{tr}(X V)) - \mathrm{Re}(\mathrm{tr}(Y \overline{W})). \label{sp trace form}
\end{align}
From this and the formula for $\n$ in \eqref{ndef}, it is also straightforward to verify that the map $\n \rightarrow \widehat{\n}$ given by sending $$\begin{pmatrix} 0 & x_1 & i\alpha_1\\ x_2 & 0 & -\ol{x_1} \\ i\alpha_2 & -\ol{x_2} & 0 \end{pmatrix} \mapsto \begin{pmatrix} 0 & x_1 & \alpha_1\\ x_2 & 0 & \ol{x_1} \\ -\alpha_2 & \ol{x_2} & 0 \end{pmatrix} \,\, \text{ and } \,\, \begin{pmatrix} y_1 & y_2 & 0\\ y_3 & 0 & y_2\\ 0 & y_3 & y_1 \end{pmatrix} \mapsto \begin{pmatrix} y_1 & y_2 & 0\\ y_3 & 0 & -y_2\\ 0 & -y_3 & -y_1 \end{pmatrix},$$ 
with $x_1,x_2,y_1,y_2,y_3 \in \C, \alpha_1,\alpha_2 \in \R$,
is an isometry.

For part (b), it suffices to show, for any $A, B \in \n$, that
\begin{align}
c^2 K_{\g}(A,B) = K_{\so(\m)}(\rho_-(A),\rho_+(B)), \label{want to show}
\end{align}
for some constant $c$, where $\rho_-(A)$ and $\rho_+(B)$ denote the projections onto the indicated grading components (i.e. onto $\p_-$, respectively $\p_+$) of $\rho(A), \rho(B) \in \so(\m)$. For, if we know that (\ref{want to show}) holds, we can simply take $\{ e_1,\ldots,e_{12}\}$ to be any $K_{\g}$-orthonormal basis of $\n$, and define $$E_i := \frac{1}{c}\rho_-(e_i) \,\, \text{ and } \,\, E^i := \frac{\varepsilon_i}{c}\rho_+(e_i).$$ Then, by construction, the relation (\ref{relate ei to Ei}) holds, and a quick calculation using (\ref{want to show}) shows that $K_{\so(\m)}(E_i,E^j) = \delta_{ij}$ for all $1 \leq i,j \leq 12$, which proves part (b).

To carry out the calculation of (\ref{want to show}), let us fix the Killing form $K_{\so(\m)}$ to be one-half the trace form over $\m$. If we write elements of $\so(\m)$ in matrix form with respect to a basis $\{ S, e_1, \ldots, e_{12}, -\widehat{S} \}$, where $\{ e_1,\ldots,e_{12} \}$ is any orthonormal basis of $\widehat{\n}$, then it is a straightforward calculation to verify that an arbitrary element of $\p_-$ has the form
\begin{align}
A_- = \begin{pmatrix} 0 & 0 & 0 \\
			x & 0 & 0 \\
			0 & -x^{\top}\1_{5,7} & 0 \end{pmatrix}. \label{p- matrices}
\end{align}
for some vector $x=(x^1,\ldots, x^{12})$ and $x=\sum_{i=1}^{12} x^ie_i \in \widehat{\n}$,  while an element of $\p_+$ has the form
\begin{align}
B_+ = \begin{pmatrix} 0 & -y^{\top}\1_{5,7} & 0 \\
			0 & 0 & y \\
			0 & 0 & 0 \end{pmatrix}. \label{p+ matrices}
\end{align}
for some vector $y \in \widehat{\n}$. In particular, we can calculate that $$K_{\so(\m)}(A_-,B_+) = -x^{\top}\mathrm{I}_{5,7}y = -K_{\g}(x,y).$$ In particular, since we have $x = A_-(S)$ and $y = -B_+(\widehat{S})$, we get the following observation: If $\wt{A}, \wt{B} \in \so(\m)$ are two elements such that $\wt A(S),\wt B(\widehat{S}) \in \widehat{\n}$, then $$K_{\so(\m)}(A_-,B_+) = K_{\g}(\wt A(S),\wt B(\widehat{S})),$$ where $A_-, B_+$ denote the projections of $\wt A$ and $\wt B$ onto the indicated grading components.

We will apply this to $\wt A =\rho(A)$, for $A \in \n$.  It is a straightforward calculation to verify that $\wt A(S)=\rho(A)S=
[A,S] \in \widehat{\n}$ and that  $\wt A (\widehat{S})=[A,\widehat{S}]\in \widehat{\n}$. Thus, the above observation applied to  $A,B\in \n$ becomes
$$K_{\so(\m)}(\rho_-(A),\rho_+(B) ) = K_{\g}([A,S],[ B,\widehat{S}]).$$
Thus, we can verify the identity (\ref{want to show}) by comparing $K_{\g}([A,S],[B,\widehat{S}])$ with $K_{\g}(A,B)$ for arbitrary $A, B \in \n$. Here are the details of that calculation:
 We let $$A = \begin{pmatrix}X&-Y\\\ol{Y}&\ol{X}\end{pmatrix}, \,\, B = \begin{pmatrix}V&-W\\\ol{W}&\ol{V}\end{pmatrix} \in \n,$$ where the matrices $X, Y, V, W \in \gl_3\C$ are given, respectively, by
\begin{align*}
X = \begin{pmatrix}
0 & x_1 & i\alpha_1\\
x_2 & 0 & -\ol{x_1} \\
i\alpha_2 & -\ol{x_2} & 0
\end{pmatrix}, \,\,
Y = \begin{pmatrix}
y_1 & y_2 & 0\\
y_3 & 0 & y_2\\
0 & y_3 & y_1
\end{pmatrix}, \,\,
V = \begin{pmatrix}
0 & v_1 & i\beta_1\\
v_2 & 0 & -\ol{v_1} \\
i\beta_2 & -\ol{v_2} & 0
\end{pmatrix}, \,\,
W = \begin{pmatrix}
w_1 & w_2 & 0\\
w_3 & 0 & w_2\\
0 & w_3 & w_1
\end{pmatrix}.
\end{align*}

Recalling the definition of $S$ via $S_0$ in \eqref{Seqn} and that $\widehat{S} = -\ol{S}$, a simple calculation shows
\begin{align}
[A,S] = \begin{pmatrix} [X,S_0] & -Y\overline{S_0} + S_0Y \\ \overline{Y}S_0 - \overline{S_0}\overline{Y} & [\overline{X},\overline{S_0}] \end{pmatrix}, \,\,
[B,\widehat{S}] = \begin{pmatrix} -[V,\ol{S_0}] & WS_0 - \ol{S_0}W \\ -\ol{W}\ol{S_0} + S_0\ol{W} & -[\ol{V},S_0] \end{pmatrix}. \label{bracket with S}
\end{align}
Furthermore, substituting $S_0=\mathrm{diag}(\mu,-2a,\overline{\mu})$ with $a=\mathrm{Re}(\mu)$, we get
\begin{align}
[X,S_0] &= \begin{pmatrix} 0 & -x_1(\mu + 2a) & i\alpha_1(\overline{\mu}-\mu) \\
						x_2(\mu + 2a) & 0 & -\overline{x_1}(\overline{\mu}+2a) \\
						i\alpha_2(\mu-\overline{\mu}) & \overline{x_2}(\overline{\mu}+2a) & 0 \end{pmatrix}; \label{[A,S] formula 1}\\
Y\overline{S_0} - S_0Y &= \begin{pmatrix} y_1(\overline{\mu}-\mu) & -y_2(\mu+2a) & 0\\
								y_3(\overline{\mu}+2a) & 0 & y_2(\mu + 2a) \\
								0 & -y_3(\overline{\mu}+2a) & -y_1(\overline{\mu}-\mu) \end{pmatrix}; \label{[A,S] formula 2}\\
-[V,\ol{S_0}] &= \begin{pmatrix} 0 & v_1(\ol{\mu}+2a) & i\beta_1(\ol{\mu}-\mu) \\
							-v_2(\ol{\mu}+2a) & 0 & \ol{v_1}(\mu+2a) \\
							i\beta_2(\mu-\ol{\mu}) & -\ol{v_2}(\mu+2a) & 0 \end{pmatrix}; \label{[B,S] formula 1}\\
-WS_0 + \ol{S_0}W &= \begin{pmatrix} w_1(\ol{\mu}-\mu) & w_2(\ol{\mu}+2a) & 0 \\
								-w_3(\mu+2a) & 0 & -w_2(\ol{\mu}+2a) \\
								0 & w_3(\mu+2a) & -w_1(\ol{\mu}-\mu) \end{pmatrix}. \label{[B,S] formula 2}
\end{align}

From the formula (\ref{sp trace form}), we therefore have
\begin{align}
K_{\g}(A,B) &= \mathrm{Re}(\mathrm{tr}(XV)) - \mathrm{Re}(\mathrm{tr}(Y\ol{W})); \label{KAB formula 1}\\
K_{\g}([A,S],[B,\widehat{S}]) &= -\mathrm{Re}(\mathrm{tr}\left([X,S_0] \circ [V,\ol{S_0}]\right))
		+ \mathrm{Re}(\mathrm{tr}\left((Y\ol{S_0} - S_0Y) \circ (\ol{W}\ol{S_0} - S_0\ol{W})\right)). \label{KASBS formula 1}
\end{align}

But now, using the form of the matrices $X, Y, V, W$ from above, we can calculate the traces in the right-hand side of (\ref{KAB formula 1}), to get:
\begin{align}
K_{\g}(A,B) &= 2\mathrm{Re}(x_1v_2 + x_2v_1 - y_1\ol{w_1} - y_2\ol{w_3} - y_3\ol{w_2}) - (\alpha_1\beta_2 + \alpha_2\beta_1). \label{KAB final formula}
\end{align}

Similarly, we can use the formulas (\ref{[A,S] formula 1}) - (\ref{[B,S] formula 2}) to compute the traces in the right-hand side of (\ref{KASBS formula 1}), to get:
\begin{align*}
K_{\g}([A,S],[B,\widehat{S}]) &= 2 \vert \mu + 2a \vert^2 \mathrm{Re}(x_1 v_2 + x_2 v_1) + (\mu-\ol{\mu})^2(\alpha_1\beta_2 + \alpha_2\beta_1) \\
&+ 2(\mu-\ol{\mu})^2\mathrm{Re}(y_1\ol{w_1}) -2\vert \mu + 2a \vert^2\mathrm{Re}(y_2\ol{w_3} + y_3\ol{w_2}).
\end{align*}
We can simplify the right-hand side of this formula by noticing that, for $\mu = a + ib$, the matrix $S$ is null with respect to $K_{\g}$ precisely when $b^2 = 3a^2$. Thus, we see that $\vert \mu + 2a \vert^2 = 9a^2 + b^2 = 12a^2$, while $(\mu - \ol{\mu})^2 = -4b^2 = -12a^2$, and hence the last display simplifies to
\begin{align}
K_{\g}([A,S],[B,\widehat{S}]) &= 24a^2\mathrm{Re}(x_1 v_2 + x_2 v_1 - y_1\ol{w_1} - y_2\ol{w_3} - y_3\ol{w_2}) -12a^2(\alpha_1\beta_2 + \alpha_2\beta_1). \label{KASBS final formula}
\end{align}
Therefore, comparing (\ref{KAB final formula}) with (\ref{KASBS final formula}), we see that $$12a^2K_{\g}(A,B) = K_{\g}([A,S],[B,\widehat{S}]),$$ as required. \end{proof}

\end{appendix}



\def\cprime{$'$}

\end{document}